\documentclass[11pt]{amsart}
\usepackage{amsmath,amssymb,latexsym,esint,cite,mathrsfs}
\usepackage{verbatim,wasysym}
\usepackage[left=2.6cm,right=2.6cm,top=2.9cm,bottom=2.9cm]{geometry}
\usepackage{xcolor}

\usepackage{graphicx}
\usepackage{subfigure}

\usepackage{graphicx,epic,eepic}

\allowdisplaybreaks

\begin{document}

\title[Fractional Hamiltonian systems]
{Existence of solution for a class of fractional Hamiltonian-type elliptic systems with exponential critical growth in $\mathbb{R}$}

\author[S. Deng]{Shengbing Deng}
\address{\noindent S. Deng-School of Mathematics and Statistics, Southwest University,
Chongqing 400715, People's Republic of China}\email{shbdeng@swu.edu.cn}

\author[J. Yu]{Junwei Yu}
\address{\noindent J. Yu-School of Mathematics and Statistics, Southwest University,
Chongqing 400715, People's Republic of China.}\email{JwYumaths@163.com}

\maketitle

\maketitle
\numberwithin{equation}{section}
\newtheorem{theorem}{Theorem}[section]
\newtheorem{lemma}[theorem]{Lemma}
\newtheorem{corollary}[theorem]{Corollary}
\newtheorem{definition}[theorem]{Definition}
\newtheorem{proposition}[theorem]{Proposition}
\newtheorem{remark}[theorem]{Remark}
\allowdisplaybreaks

\maketitle

\noindent {\bf Abstract}: In this paper, we study the following class of fractional Hamiltonian systems:
\begin{eqnarray*}
    \begin{aligned}\displaystyle
    \left\{ \arraycolsep=1.5pt
       \begin{array}{ll}
       (-\Delta)^{\frac{1}{2}} u + u = \Big(I_{\mu_{1}}\ast G(v)\Big)g(v) \ \ \ & \mbox{in} \ \mathbb{R},\\[2mm]
       (-\Delta)^{\frac{1}{2}} v + v = \Big(I_{\mu_{2}}\ast F(u)\Big)f(u) \ \ \ & \mbox{in} \ \mathbb{R},
        \end{array}
    \right.
    \end{aligned}
\end{eqnarray*}
where $(-\Delta)^{\frac{1}{2}}$ is the square root Laplacian operator, $\mu_{1},\mu_{2}\in(0,1)$, $I_{\mu_{1}},I_{\mu_{2}}$ denote the Riesz potential, $\ast$ indicates the convolution operator, $F(s),G(s)$ are the primitive of $f(s),g(s)$ with $f(s),g(s)$ have exponential growth in $\mathbb{R}$. Using the linking theorem and variational methods, we establish the existence of at least one positive solution to the above problem.

\vspace{3mm} \noindent {\bf Keywords}: Fractional Hamiltonian elliptic system; Choquard nonlinearities; Critical exponential growth; Trudinger-Moser inequality.

\vspace{3mm}

\maketitle

\section{{\bfseries Introduction}}

This paper deals with the following class of fractional Hamiltonian elliptic systems of the Choquard type
\begin{equation}\label{a}
    \begin{aligned}\displaystyle
    \left\{ \arraycolsep=1.5pt
       \begin{array}{ll}
       (-\Delta)^{\frac{1}{2}} u + u = \Big(I_{\mu_{1}}\ast G(v)\Big)g(v) \ \ \ & \mbox{in} \ \mathbb{R},\\[2mm]
       (-\Delta)^{\frac{1}{2}} v + v = \Big(I_{\mu_{2}}\ast F(u)\Big)f(u) \ \ \ & \mbox{in} \ \mathbb{R},
        \end{array}
    \right.
    \end{aligned}
\end{equation}
where $(-\Delta)^{\frac{1}{2}}$ denotes the fractional Laplacian, $\mu_{1},\mu_{2}\in(0,1)$, $I_{\mu_{1}},I_{\mu_{2}}$ denote the Riesz potential defined by
\begin{eqnarray*}
    \begin{aligned}\displaystyle
    I_{\mu_{i}}(x)=\frac{\Gamma(\frac{2-\mu_{i}}{2})}{\Gamma(\frac{\mu_{i}}{2})2^{\mu_{i}}\pi|x|^{2-\mu_{i}}}:=\frac{A_{\mu_{i}}}{|x|^{2-\mu_{i}}}, \ \ x\in \mathbb{R}\backslash\{0\}, \ i=1,2,
    \end{aligned}
\end{eqnarray*}
where $\Gamma$ represents the gamma function, $\ast$ indicates the convolution operator, $F(t)=\int^{t}_{0}f(\tau)d\tau$ and $G(t)=\int^{t}_{0}g(\tau)d\tau$. The nonlinearities $f,g$ satisfy some suitable conditions that will be specified later. The operator $(-\Delta)^{\frac{1}{2}}$ is the square root Laplacian operator defined as
\begin{eqnarray*}
    \begin{aligned}\displaystyle
    (-\Delta)^{\frac{1}{2}} u(x)=-\frac{1}{2\pi}\int^{1}_{-1}\frac{u(x+y)+u(x-y)-2u(x)}{|y|^{2}}dy,
    \end{aligned}
\end{eqnarray*}

We say that a function $h(t)$ has critical exponential growth if there exists $\gamma_{0}>0$ such that
\begin{equation}\label{ac}
    \begin{aligned}\displaystyle
    \lim\limits_{|t|\rightarrow +\infty}\frac{|h(t)|}{e^{\gamma s^{2}}}=0,\ \ \forall \ \gamma>\gamma_{0},\ \ and \ \ \lim\limits_{|t|\rightarrow +\infty}\frac{|h(t)|}{e^{\gamma t^{2}}}=+ \infty,\ \ \forall \ \gamma<\gamma_{0}.
    \end{aligned}
\end{equation}
This notion of criticality is motivated by a class of Trudinger-Moser inequality, see \cite{Ta,Lula}.

In recent years, there are many works dedicated to study the following non-local elliptic equation
\begin{equation}\label{b}
    \begin{aligned}\displaystyle
    (-\Delta)^{s} u +u =(I_{\mu} \ast F(u))f(u) \ \ \mbox{in} \ \mathbb{R}^{N},
    \end{aligned}
\end{equation}
where $(-\Delta)^{s}$ denotes the  fractional Laplacian with $0<s<1$,  $I_{\mu}=\frac{1}{|x|^{\mu}}$ denotes the Riesz potential. This nonlocal equation arises in many interesting physical situations in quantum theory and plays an important role in describing the finite-range many-body interactions. When $s=1$ and $F(t)=\frac{1}{2}t^{2}$, equation $(\ref{b})$ was introduced by Fr\"{o}hlich in \cite{FR} to study the modeling of quantum polaron. We mention that a great attention has been focused on the study of problem involving exponential growth nonlinearity. We know that classical fractional Sobolev embedding that $H^{s,2}(\mathbb{R}^{N})$ is continuously embedded in $L^{q}(\mathbb{R}^{N})$ for all $q\in [2,2_{s}^{*}]$, where $2_{s}^{*}=2N/(N-2s)$. Thus, we know that $2_{s}^{*}=\infty$ if $N=2s$. In the case $N=1$ and $s=1/2$, Clemente et al. \cite{Cle} studied the existence of solutions of the following equation
\begin{equation}\label{ba}
    \begin{aligned}\displaystyle
    (-\Delta)^{\frac{1}{2}} u +u =(I_{\mu} \ast F(u))f(u) \ \ \mbox{in} \ \mathbb{R},
    \end{aligned}
\end{equation}
where $f$ has critical exponential growth in the sense of Trudinger-Moser inequality.

In the case $s=1$, Alves et al. \cite{ACTT} obtained the existence and concentration results for a class of nonlinear Choquard equations in the plane with exponential critical nonlinearity. Battaglia and Van Schaftingen \cite{BATTA} proved the existence of a nontrivial ground state solution for equation $(\ref{b})$ when $f$ has exponential growth. x Qin and Tang \cite{QT} improved and extended the related results to the strongly indefinite problems. For more classical results regarding Choquard equations, the reader may refer \cite{MVb} to the good survey.

On the other hand, we mention that the study of our problem is based on some interesting results of fractional Hamiltonian elliptic systems.  do \'{O} et al. \cite{Do} considered the following nonautonomous fractional Hamiltonian system with critical exponential growth
\begin{equation}\label{c}
    \begin{aligned}\displaystyle
    \left\{ \arraycolsep=1.5pt
       \begin{array}{ll}
       (-\Delta)^{\frac{1}{2}} u + u = Q(x)g(v) \ \ \ & \mbox{in} \ \mathbb{R},\\[2mm]
       (-\Delta)^{\frac{1}{2}} v + v = P(x)f(u) \ \ \ & \mbox{in} \ \mathbb{R},
        \end{array}
    \right.
    \end{aligned}
\end{equation}
where $f$ and $g$ have critical growth at infinity in the sense of Trudinger-Moser inequality and the nonnegative weights $P(x)$ and $Q(x)$ vanish at infinity. Combining Moser functions and the asymptotic behavior of $f,g$ at infinity, the authors of \cite{Do} estimated the minimax level of functional associated with $(\ref{c})$, by using suitable variational method combined with the generalized linking theorem, they obtained the existence of at least one positive solution for $(\ref{c})$.

When $s=1$ and $N=2$, a great attention has been focused on the study of the following Hamiltonian elliptic system
\begin{equation}\label{cc}
    \begin{aligned}\displaystyle
    \left\{ \arraycolsep=1.5pt
       \begin{array}{ll}
       -\Delta u + V(x)u = f_{1}(x,v) \ \ \ & \mbox{in} \ \Omega,\\[2mm]
       -\Delta v + V(x)v = f_{2}(x,u) \ \ \ & \mbox{in} \ \Omega,
        \end{array}
    \right.
    \end{aligned}
\end{equation}
where $\Omega\subset \mathbb{R}^{2}$ is a bound domain. For $V=0$, de Figueiredo et al. \cite{FdR} studied the existence of solution for $(\ref{cc})$ with the Ambrosetti-Rabinowitz condition. Lam and Lu \cite{LLU} considered the existence of nontrivial nonnegative solutions for $(\ref{cc})$ with subcritical and critical exponential growth without the Ambrosetti-Rabinowitz condition. When $V=1$ and $\Omega=\mathbb{R}^{2}$, the existence of ground state solutions of $(\ref{cc})$ was introduced by de Figueiredo et al. \cite{FJZ}. For more classical results regarding Hamiltonian system, we refer to \cite{Bo} a good survey.

When the nonlinear term is Choquard type, Maia and Miyagaki \cite{MMA} considered the following system
\begin{equation}\label{e}
    \begin{aligned}\displaystyle
    \left\{ \arraycolsep=1.5pt
       \begin{array}{ll}
       -\Delta u + V(x)u = \Big(I_{\mu}\ast G(v)\Big)g(v) \ \ \ & \mbox{in} \ \mathbb{R}^{2},\\[2mm]
       -\Delta v + V(x)v = \Big(I_{\mu}\ast F(v)\Big)f(v) \ \ \ & \mbox{in} \ \mathbb{R}^{2},
        \end{array}
    \right.
    \end{aligned}
\end{equation}
where the potential $V$ and the nonlinearities $f,g$ satisfy some suitable conditions. Using approximation methods and linking arguments, they showed the existence of solutions for $(\ref{e})$ with critical exponential growth.
Recently, Tang et al. \cite{QTZA} obtained the existence of ground states of $(\ref{cc})$ by the non-Nehari manifold method. They developed a direct method to deal with the difficulties aroused by the strongly indefinite features and the critical exponential growth.

The purpose of this work is to consider the existence of solutions for the fractional Hamiltonian Choquard-type elliptic systems $(\ref{a})$. Assume that $f(t)=g(t)=0$ for $t\leq0$, and the nonlinearities $f,g$ satisfy the following conditions:

$(H_{0})$ $f,g: \mathbb{R} \rightarrow [0,+\infty)$ are continuous functions;

$(H_{1})$ $\lim\limits_{t\rightarrow 0^{+}}f(t)t^{-1}=\lim\limits_{s\rightarrow 0^{+}}g(t)t^{-1}=0$;

$(H_{2})$ there exists $\theta > 1$ such that for all $t>0$,
    \begin{eqnarray*}
    \begin{aligned}\displaystyle
    0<\theta F(t)\leq tf(t) \ \ and \ \ 0< \theta G(t) \leq tg(t);
    \end{aligned}
    \end{eqnarray*}

$(H_{3})$ there exist constants $t_{0},M_{0}>0$ such that for all $t>t_{0}$,
    \begin{eqnarray*}
    \begin{aligned}\displaystyle
    0<F(t)\leq M_{0}f(t) \ \ and \ \ 0< G(t) \leq M_{0}g(t);
    \end{aligned}
    \end{eqnarray*}

$(H_{4})$ There exist $\alpha_{0}>0$ and $\beta_{0}>0$ such that
\begin{eqnarray*}
    \begin{aligned}\displaystyle
    \lim\limits_{|t|\rightarrow\infty}\frac{|f(t)|}{e^{\alpha t^{2}}}=
    \left\{ \arraycolsep=1.5pt
       \begin{array}{ll}
       +\infty, \ \ \ & \alpha<\alpha_{0},\\[2mm]
       0, \ \ \ & \alpha>\alpha_{0},
        \end{array}
    \right.
    \ \ \ \mbox{and}\ \ \     \lim\limits_{|t|\rightarrow\infty}\frac{|g(t)|}{e^{\beta t^{2}}}=
    \left\{ \arraycolsep=1.5pt
       \begin{array}{ll}
       +\infty, \ \ \ & \beta<\beta_{0},\\[2mm]
       0, \ \ \ & \beta>\beta_{0};
        \end{array}
    \right.
    \end{aligned}
\end{eqnarray*}

$(H_{5})$ $\displaystyle\liminf\limits_{t\rightarrow +\infty}\frac{tF(t)}{e^{\alpha_{0}t^{2}}}\geq \kappa_{1}$ with $\displaystyle\kappa_{1}>\sqrt{\frac{\pi\mu_{2}(1+\mu_{2})^{2}}{2^{3+\mu_{2}}A_{\mu_{2}}\alpha_{0}^{2}}e^{\frac{4(1+\mu_{2})}{\pi}-1}}$,

\quad \quad \ and \quad $\displaystyle\liminf\limits_{t\rightarrow +\infty}\frac{tG(t)}{e^{\beta_{0}t^{2}}}\geq \kappa_{2}$ with $\displaystyle\kappa_{2}> \sqrt{\frac{\pi\mu_{1}(1+\mu_{1})^{2}}{2^{3+\mu_{1}}A_{\mu_{1}}\beta_{0}^{2}}e^{\frac{4(1+\mu_{1})}{\pi}-1}}$.

The main result can be statedx as follows.

\begin{theorem}\label{THb}
Suppose that $(H_{0})-(H_{5})$ hold, then system $(\ref{a})$ has a positive weak solution.
\end{theorem}

\begin{remark}\rm
There are some difficulties in system $(\ref{a})$. The first one is that the energy functional associated with system $(\ref{a})$ is strongly indefinite and there is a lack of the compactness for the whole space case. Hence, we use an approximation procedure to overcome the lack of compactness. The second one is the nonlinearities with critical exponents $\alpha_{0}$ and $\beta_{0}$. Many difficulties arise when $\alpha_{0}$ is not equal to $\beta_{0}$. In order to address this problem, we use the definitions of hole functions in \cite{Ca,Milan}. Using the Moser type functions and the hole functions, we restore compactness by the properties of the exponential critical nonlinear term at infinity.
\end{remark}
The paper is organized as follows: In Section $\ref{PFS}$, the variational setting and some preliminary basic results are presented. In Section $\ref{STLEMMA}$, we introduce the geometry structure and give some properties related to system $(\ref{a})$. In Section $\ref{FIn}$, we give an approximation procedure of system $(\ref{a})$. Section $\ref{TECL}$ is devoted to give an estimate for the minimax level of the critical case. In Section $\ref{PROOF}$, we prove Theorem $\ref{THb}$.

\section{{\bfseries Preliminaries and functional setting}}\label{PFS}
In this section, we give some preliminary results and outline the variational framework for $(\ref{a})$.
\begin{proposition}\label{PRa}(A fractional Trudinger-Moser inequality)\cite{Lula}.
It holds
\begin{eqnarray*}
    \begin{aligned}\displaystyle
    \sup\limits_{u\in H^{1/2}(\mathbb{R}),\|u\|_{1/2}\leq1}\int_{\mathbb{R}}(e^{\gamma|u|^{2}}-1)dx
    \left\{ \arraycolsep=1.5pt
       \begin{array}{ll}
       <\infty, \ \ \ & \gamma \ \leq\pi,\\[2mm]
       =\infty, \ \ \ & \gamma \ >\pi,
        \end{array}
    \right.
    \end{aligned}
\end{eqnarray*}
where $H^{1/2}(\mathbb{R})$ is the fractional order Sobolev space equipped with $\|\cdot\|_{1/2}$ norm which is defined in the following page.
\end{proposition}

Since we are going to study the nonlocal type problems with the Riesz potential, we would like to state Hardy-Littlewood-Sobolev inequality in $\mathbb{R}$.

\begin{proposition}\label{PRc}(Hardy-Littlewood-Sobolev inequality)\cite{Lieb}. Let $t,r>1$, $0<\mu<1$, with $\frac{1}{t}+\mu+\frac{1}{r}=2$, $f\in L^{t}(\mathbb{R})$ and $h\in L^{r}(\mathbb{R})$. Then there exists a constant $C(t,\mu,r)$, independent of $f,h$, such that
\begin{eqnarray*}
    \begin{aligned}\displaystyle
    \int_{\mathbb{R}}\int_{\mathbb{R}}\frac{f(x)h(y)}{|x-y|^{\mu}}dxdy \leq C(t,\mu,r)\|f\|_{L^{t}(\mathbb{R})}\|h\|_{L^{r}(\mathbb{R})}.
    \end{aligned}
\end{eqnarray*}
\end{proposition}

Applying the Hardy-Littlewood-Sobolev inequality, we know
\begin{eqnarray*}
    \begin{aligned}\displaystyle
    \int_{\mathbb{R}}\Big(I_{\mu_{2}}\ast F(u)\Big)F(u)dx, \ \ \int_{\mathbb{R}}\Big(I_{\mu_{1}}\ast G(u)\Big)G(u)dx
    \end{aligned}
\end{eqnarray*}
are well defined if $F(u),G(u)\in L^{t}(\mathbb{R})$ for $t>1$ given by
\begin{eqnarray*}
    \begin{aligned}\displaystyle
    \frac{2}{t}+1-\mu_{2}=2, \ \ \frac{2}{t}+1-\mu_{1}=2.
    \end{aligned}
\end{eqnarray*}
This implies that we must require
\begin{eqnarray*}
    \begin{aligned}\displaystyle
    F(u)\in L^{\frac{2}{1+\mu_{2}}}(\mathbb{R}), \ \ G(u)\in L^{\frac{2}{1+\mu_{1}}}(\mathbb{R}).
    \end{aligned}
\end{eqnarray*}
In order to apply variational methods, we consider the following subspace of $H^{1/2}(\mathbb{R})$
\begin{eqnarray*}
    \begin{aligned}\displaystyle
    H^{1/2}(\mathbb{R})=\Big\{u\in L^{2}(\mathbb{R}): [u]_{1/2}<\infty\},
    \end{aligned}
\end{eqnarray*}
where the term
\begin{eqnarray*}
    \begin{aligned}\displaystyle
    [u]_{1/2}=\Big(\int_{\mathbb{R}}\int_{\mathbb{R}}\frac{|u(x)-u(y)|^{2}}{|x-y|^{2}}dxdy\Big)^{\frac{1}{2}},
    \end{aligned}
\end{eqnarray*}
is the so-called Gagliardo semi-norm of a function $u$. We point out from \cite[Proposition 3.6]{Di} that
\begin{eqnarray*}
    \begin{aligned}\displaystyle
    \|(-\Delta)^{1/4}u\|^{2}_{L^{2}(\mathbb{R})}=\frac{1}{2\pi}\int_{\mathbb{R}}\int_{\mathbb{R}}\frac{|u(x)-u(y)|^{2}}{|x-y|^{2}}dxdy, \ \ \ \mbox{for} \ \mbox{all} \ \ u\in H^{1/2}(\mathbb{R}).
    \end{aligned}
\end{eqnarray*}
The space $H^{1/2}(\mathbb{R})$ is a Hilbert space when endowed with the inner product
\begin{eqnarray*}
    \begin{aligned}\displaystyle
     \langle u,v \rangle_{1/2}=\int_{\mathbb{R}}[(-\Delta)^{\frac{1}{4}} u(-\Delta)^{\frac{1}{4}} v + uv]dx, \ \ u,v\in H^{1/2}(\mathbb{R})
    \end{aligned}
\end{eqnarray*}
which corresponds the norm
\begin{eqnarray*}
    \begin{aligned}\displaystyle
    \|u\|_{1/2}=\Big(\frac{1}{2\pi}[u]_{1/2}^{2}+\int_{\mathbb{R}} u^{2}dx\Big)^{\frac{1}{2}}.
    \end{aligned}
\end{eqnarray*}

Set $E:=H^{1/2}(\mathbb{R})\times H^{1/2}(\mathbb{R})$. Then E endowed with the inner product
\begin{eqnarray*}
    \begin{aligned}\displaystyle
     \langle (u,v),(\varphi,\psi) \rangle_{E}=\int_{\mathbb{R}}[(-\Delta)^{\frac{1}{4}} u(-\Delta)^{\frac{1}{4}} \varphi +(-\Delta)^{\frac{1}{4}} v(-\Delta)^{\frac{1}{4}} \psi + u\varphi +v\psi]dx
    \end{aligned}
\end{eqnarray*}
for all $(u,v),(\varphi,\psi)\in E$, to which corresponds the norm $\|(u,v)\|_{E}=\langle (u,v), (u,v)\rangle_{E}^{\frac{1}{2}}$.

The energy functional associated with $(\ref{a})$ is given by $\Phi:E\rightarrow \mathbb{R}$, where
\begin{eqnarray*}
    \begin{aligned}\displaystyle
    \Phi(u,v)=\langle u,v\rangle_{1/2}-\frac{1}{2}\int_{\mathbb{R}}\Big(I_{\mu_{2}}\ast F(u)\Big)F(u)dx-\frac{1}{2}\int_{\mathbb{R}}\Big(I_{\mu_{1}}\ast G(v)\Big)G(v)dx.
    \end{aligned}
\end{eqnarray*}
Using standard arguments it is possible to verify that $\Phi$ is well defined and is of class $C^{1}$ with
\begin{eqnarray*}
    \begin{aligned}\displaystyle
    \Phi^{\prime}(u,v)(\varphi,\psi)=&\int_{\mathbb{R}}[(-\Delta)^{\frac{1}{4}} u(-\Delta)^{\frac{1}{4}} \psi+(-\Delta)^{\frac{1}{4}} v(-\Delta)^{\frac{1}{4}} \varphi + u\psi +v\varphi]dx\\
    &-\int_{\mathbb{R}}\Big(I_{\mu_{2}}\ast F(u)\Big)f(u)\varphi dx-\int_{\mathbb{R}}\Big(I_{\mu_{1}}\ast G(v)\Big)g(v)\psi dx,
    \end{aligned}
\end{eqnarray*}
for all $(\varphi,\psi)\in E$.

Let us define the following subspaces of $E$,
\begin{eqnarray*}
    \begin{aligned}\displaystyle
    E^{+}=\{(u,u)\in E\}\ \ and\ \  E^{-}=\{(v,-v)\in E\}
    \end{aligned}
\end{eqnarray*}
with
\begin{eqnarray*}
    \begin{aligned}\displaystyle
    (u,v)=\frac{1}{2}(u+v,u+v)+\frac{1}{2}(u-v,v-u),
    \end{aligned}
\end{eqnarray*}
then $E^{+}$ is orthogonal to $E^{-}$ with respect to the inner product $\langle \cdot,\cdot\rangle$. For every $z=(u,v)$, set
\begin{eqnarray*}
    \begin{aligned}\displaystyle
    z^{+}:=\Big(\frac{u+v}{2},\frac{u+v}{2}\Big) \ and \ z^{-}:=\Big(\frac{u-v}{2},\frac{v-u}{2}\Big).
    \end{aligned}
\end{eqnarray*}
Note that $z^{+}\in E^{+}$ and $z^{-}\in E^{-}$, then we have $E=E^{+}\bigoplus E^{-}$.

From $(H_{0}),(H_{1})$ and $(H_{4})$, we have the following result: for all $\alpha>\alpha_{0}$, $\beta>\beta_{0}$ and $q\geq1$, for any given $\varepsilon >0$, there exist constants $b_{1},b_{2},b_{3},b_{4}>0$ such that
\begin{equation}\label{fnon}
    \begin{aligned}\displaystyle
    f(t)\leq\varepsilon|t|+b_{1}|t|^{q-1}(e^{\alpha t^{2}}-1),\ \ \mbox{for} \ \mbox{all} \ t\in \mathbb{R},\\
    g(t)\leq\varepsilon|t|+b_{2}|t|^{q-1}(e^{\beta t^{2}}-1),\ \ \mbox{for} \ \mbox{all} \ t\in \mathbb{R},
    \end{aligned}
\end{equation}
and
\begin{equation}\label{fnonb}
    \begin{aligned}\displaystyle
    F(t)\leq\frac{\varepsilon}{2}|t|^{2}+b_{3}|t|^{q}(e^{\alpha t^{2}}-1),\ \ \mbox{for} \ \mbox{all} \ t\in \mathbb{R},\\
    G(t)\leq\frac{\varepsilon}{2}|t|^{2}+b_{4}|t|^{q}(e^{\beta t^{2}}-1),\ \ \mbox{for} \ \mbox{all} \ t\in \mathbb{R}.
    \end{aligned}
\end{equation}

The vanishing lemma was proved originally by P.L. Lions \cite[Lemma I.1]{Lions} and there we use the following version to fractional Sobolev spaces:
\begin{lemma}\label{LEbd}Assume that $(u_{n})$ is a bounded sequence in $H^{1/2}(\mathbb{R})$ satisfying
\begin{eqnarray*}
    \begin{aligned}\displaystyle
    \lim\limits_{n\rightarrow+\infty}\sup\limits_{y\in \mathbb{R}} \int^{y+R}_{y-R}|u_{n}|^{2}dx =0,
    \end{aligned}
\end{eqnarray*}
for some $R>0$. Then, $u_{n}\rightarrow 0$ strongly in $L^{p}(\mathbb{R})$, for $2<p<\infty$.
\end{lemma}

\section{{\bfseries On the Palais-Smale sequence}}\label{STLEMMA}
\begin{lemma}\label{LEca}The following inequality holds:
\begin{eqnarray*}
    \begin{aligned}\displaystyle
    st\leq
    \left\{ \arraycolsep=1.5pt
       \begin{array}{ll}
        (e^{t^{2}}-1)+s(\log s)^{\frac{1}{2}}, \ \ \ & for \ all \ t\geq0 \ and \ s\geq e^{\frac{1}{4}};\\[2mm]
        (e^{t^{2}}-1)+\frac{1}{2}s^{2} \ \ \ & for \ all \ t\geq0 \ and \ 0\leq s\leq e^{\frac{1}{4}}.\\[2mm]
        \end{array}
    \right.
    \end{aligned}
\end{eqnarray*}
\end{lemma}
\begin{proof}
\cite[Lemma 2.4]{FJZ}.
\end{proof}

\begin{lemma}\label{LEcb} Assume that $(H_{0})-(H_{2})$ and $(H_{4})$ hold. If $(u_{n},v_{n})\in E$ such that
\begin{equation}\label{PSa}
    \begin{aligned}\displaystyle
    \Phi(u_{n},v_{n})\rightarrow c \ \ and \ \ \Phi^{\prime}(u_{n},v_{n})(\varphi,\psi)\rightarrow 0 \ \ for \ all \ (\varphi,\psi)\in E,
    \end{aligned}
\end{equation}
then $(u_{n},v_{n})$ is bounded in $E$, and there exists a constant $C>0$ such that
\begin{equation}\label{PSb}
    \begin{aligned}\displaystyle
    \int_{\mathbb{R}}\Big(I_{\mu_{2}}\ast F(u_{n})\Big)f(u_{n})u_{n}dx\leq C,\ \ \ \ \int_{\mathbb{R}}\Big(I_{\mu_{1}}\ast G(v_{n})\Big)G(v_{n})v_{n}dx\leq C,
    \end{aligned}
\end{equation}
\begin{equation}\label{PSc}
    \begin{aligned}\displaystyle
    \int_{\mathbb{R}}\Big(I_{\mu_{2}}\ast F(u_{n})\Big)F(u_{n})dx\leq C,\ \ \ \ \int_{\mathbb{R}}\Big(I_{\mu_{1}}\ast G(v_{n})\Big)G(v_{n})dx\leq C.
    \end{aligned}
\end{equation}
Moreover, for every $\phi \in C^{\infty}_{0}(\mathbb{R})$,
\begin{equation}\label{PSd}
    \begin{aligned}\displaystyle
    &\lim\limits_{n\rightarrow\infty}\int_{\mathbb{R}}\Big(I_{\mu_{2}}\ast F(u_{n})\Big)f(u_{n})\phi dx=\int_{\mathbb{R}}\Big(I_{\mu_{2}}\ast F(u)\Big)f(u)\phi dx,\\
    &\lim\limits_{n\rightarrow\infty}\int_{\mathbb{R}}\Big(I_{\mu_{1}}\ast G(v_{n})\Big)g(v_{n})\phi dx=\int_{\mathbb{R}}\Big(I_{\mu_{1}}\ast G(v)\Big)g(v)\phi dx.
    \end{aligned}
\end{equation}
\end{lemma}

\begin{proof} Taking $(\varphi,\psi)=(u_{n},v_{n})$ in $(\ref{PSa})$, we have
\begin{eqnarray*}
    \begin{aligned}\displaystyle
    2 \| (u_{n},v_{n}) \|^{2}_{E} -\int_{\mathbb{R}}\Big(I_{\mu_{2}}\ast F(u_{n})\Big)f(u_{n})u_{n}dx- \int_{\mathbb{R}}\Big(I_{\mu_{1}}\ast G(v_{n})\Big)g(v_{n})v_{n}dx = \varepsilon_{n},
    \end{aligned}
\end{eqnarray*}
which together with $(\ref{PSa})$ and hypothesis $(H_{2})$, implies
\begin{eqnarray*}
    \begin{aligned}\displaystyle
    &\int_{\mathbb{R}}\Big(I_{\mu_{2}}\ast F(u_{n})\Big)f(u_{n})u_{n}dx+ \int_{\mathbb{R}}\Big(I_{\mu_{1}}\ast G(v_{n})\Big)g(v_{n})v_{n}dx\\
    &= 2c +\varepsilon_{n} + \int_{\mathbb{R}}\Big(I_{\mu_{2}}\ast F(u_{n})\Big)F(u_{n})dx+ \int_{\mathbb{R}}\Big(I_{\mu_{1}}\ast G(v_{n})\Big)G(v_{n})dx\\
    &\leq 2c +\varepsilon_{n} +\frac{1}{\theta} \int_{\mathbb{R}}\Big(I_{\mu_{2}}\ast F(u_{n})\Big)f(u_{n})u_{n}dx+ \frac{1}{\theta} \int_{\mathbb{R}}\Big(I_{\mu_{1}}\ast G(v_{n})\Big)g(v_{n})v_{n}dx.
    \end{aligned}
\end{eqnarray*}
Thus, we obtain
\begin{eqnarray*}
    \begin{aligned}\displaystyle
    \int_{\mathbb{R}}\Big(I_{\mu_{2}}\ast F(u_{n})\Big)f(u_{n})u_{n}dx+ \int_{\mathbb{R}}\Big(I_{\mu_{1}}\ast G(v_{n})\Big)g(v_{n})v_{n}dx \leq \frac{\theta}{\theta -1} (2c +\varepsilon_{n})
    \end{aligned}
\end{eqnarray*}
where $\varepsilon_{n}\rightarrow0$, as $n\rightarrow\infty$. Next taking $(\varphi,\psi)=(v_{n},0)$ and $(\varphi,\psi)=(0,u_{n})$ in $(\ref{PSa})$ we have
\begin{eqnarray*}
    \begin{aligned}\displaystyle
    \|v_{n}\|^{2}_{1/2}\leq \int_{\mathbb{R}}\Big(I_{\mu_{2}}\ast F(u_{n})\Big)f(u_{n})v_{n}dx,
   \quad\mbox{and}\ \
    \|u_{n}\|^{2}_{1/2}\leq \int_{\mathbb{R}}\Big(I_{\mu_{1}}\ast G(v_{n})\Big)g(v_{n})u_{n}dx.
    \end{aligned}
\end{eqnarray*}
Setting $U_{n}=u_{n}/\|u_{n}\|_{1/2}$ and $V_{n}=v_{n}/\|v_{n}\|_{1/2}$, we conclude that
\begin{eqnarray*}
    \begin{aligned}\displaystyle
    \|v_{n}\|_{1/2} \leq  \int_{\mathbb{R}}\Big(I_{\mu_{2}}\ast F(u_{n})\Big)f(u_{n})V_{n}dx,
   \quad\mbox{and}\ \
    \|u_{n}\|_{1/2}\leq  \int_{\mathbb{R}}\Big(I_{\mu_{1}}\ast G(v_{n})\Big)g(v_{n})U_{n}dx.
    \end{aligned}
\end{eqnarray*}
From $(H_{1})$ and $(H_{4})$, we have
\begin{equation}\label{PSe}
    \begin{aligned}\displaystyle
    f(t)\leq C_{1}e^{\alpha t^{2}}, \ \mbox{for} \ \mbox{all} \ t\geq 0.
    \end{aligned}
\end{equation}

Using Lemma $\ref{LEca}$ with $t=V_{n}$ and $s=f(u_{n})/C_{1}$, where $C_{1}$ is the constant appearing in $(\ref{PSe})$, and the Trudinger-Moser inequality, we obtain
\begin{eqnarray*}
    \begin{aligned}\displaystyle
    C_{1}& \int_{\mathbb{R}}\Big(I_{\mu_{2}}\ast F(u_{n})\Big)\frac{f(u_{n})}{C_{1}}V_{n}dx\leq C_{1} \int_{\mathbb{R}}\Big(I_{\mu_{2}}\ast F(u_{n})\Big)(e^{V_{n}^{2}}-1)dx\\
    &+ C_{1} \int_{\{x\in \mathbb{R}:f(u_{n})/C_{1}\geq e^{1/4}\}}\Big(I_{\mu_{2}}\ast F(u_{n})\Big)\frac{f(u_{n})}{C_{1}}\Big[\log \frac {f (u_{n})}{C_{1}}\Big]^{\frac{1}{2}}dx\\
    &+ \frac{1}{2} \int_{\{x\in \mathbb{R}:f(u_{n})/C_{1}\leq e^{1/4}\}}\Big(I_{\mu_{2}}\ast F(u_{n})\Big)\frac{[f(u_{n})]^{2}}{C_{1}^{2}}dx\\
    :=& I_{1}+ I_{2} +I_{3}.
    \end{aligned}
\end{eqnarray*}
First, by the Hardy-Littlewood-Sobolev inequality, $V_{n}^{2}< \pi$, and $F(u_{n})\in L^{\frac{2}{1+\mu_{2}}}$, we obtian
\begin{eqnarray*}
    \begin{aligned}\displaystyle
    I_{1}\leq |F(u_{n})|_{\frac{2}{1+\mu_{2}}}|e^{V_{n}^{2}}-1|_{\frac{2}{1+\mu_{2}}}
    \leq C_{2}.
    \end{aligned}
\end{eqnarray*}
Next by $(\ref{PSe})$, we get
\begin{eqnarray*}
    \begin{aligned}\displaystyle
I_{2} \leq C_{3} \int_{\mathbb{R}}\Big(I_{\mu_{2}}\ast F(u_{n})\Big)f(u_{n})u_{n}dx.
    \end{aligned}
\end{eqnarray*}
Finally, using $(H_{1})$, we have
\begin{eqnarray*}
    \begin{aligned}\displaystyle
    [f(t)]^{2} \ \leq \ f(t)t,
    \end{aligned}
\end{eqnarray*}
where $t\in \{t\in \mathbb{R}:s\geq0 \ and \ f(t)/C_{1} \leq e^{1/4} \}$. Thus
\begin{eqnarray*}
    \begin{aligned}\displaystyle
    I_{3}\leq C_{4} \int_{\mathbb{R}}\Big(I_{\mu_{2}}\ast F(u_{n})\Big)f(u_{n})u_{n}dx.
    \end{aligned}
\end{eqnarray*}
Thus, we obtain that
\begin{eqnarray*}
    \begin{aligned}\displaystyle
    C_{1} \int_{\mathbb{R}}\Big(I_{\mu_{2}}\ast F(u_{n})\Big)\frac{f(u_{n})}{C_{1}}V_{n}dx\leq C_{2}+ C_{5} \int_{\mathbb{R}}\Big(I_{\mu_{2}}\ast F(u_{n})\Big)f(u_{n})u_{n}dx,
    \end{aligned}
\end{eqnarray*}
which implies that
\begin{equation}\label{PSf}
    \begin{aligned}\displaystyle
    \|v_{n}\|_{1/2}\leq C_{2} + C \int_{\mathbb{R}}\Big(I_{\mu_{2}}\ast F(u_{n})\Big)f(u_{n})u_{n}dx.
    \end{aligned}
\end{equation}
Similarly, we have
\begin{equation}\label{PSg}
    \begin{aligned}\displaystyle
    \|u_{n}\|_{1/2}\leq C_{2} + C \int_{\mathbb{R}}\Big(I_{\mu_{1}}\ast G(v_{n})\Big)g(v_{n})v_{n}dx.
    \end{aligned}
\end{equation}

Now joining the estimates $(\ref{PSf})$ and $(\ref{PSg})$ we finally obtain
\begin{eqnarray*}
    \begin{aligned}\displaystyle
    \|u_{n}\|_{1/2}+\|v_{n}\|_{1/2}\leq\frac{\theta}{\theta-1}(2c+\varepsilon_{n}),
    \end{aligned}
\end{eqnarray*}
which implies that the boundedness of $(u_{n},v_{n})$. From this estimate, we obtain $(\ref{PSb})$ and $(\ref{PSc})$. Then, by \cite[Lemma 2.1]{DDRU}, we get $(\ref{PSd})$.
\end{proof}

To find a solution of $(\ref{a})$, we show that the functional $\Phi$ possesses the following geometry.

\begin{lemma}\label{LEcd}Suppose that $(H_{0})-(H_{1})$ and $(H_{4})$ hold. Then there exists $\rho,\sigma>0$ such that
\begin{eqnarray*}
    \begin{aligned}\displaystyle
    \kappa:=\inf\{\Phi(z):z\in E^{+},\|z\|=\rho\}\geq\sigma>0.
    \end{aligned}
\end{eqnarray*}
\end{lemma}

\begin{proof} Taking $q=3$ in $(\ref{fnonb})$ and $z\in E^{+}$, we know that
\begin{eqnarray*}
    \begin{aligned}\displaystyle
    F(u)\leq\frac{\varepsilon}{2}|u|^{2}+b_{3}|u|^{3}(e^{\alpha u^{2}}-1),\ \ \mbox{and} \ \
    G(u)\leq\frac{\varepsilon}{2}|u|^{2}+b_{4}|u|^{3}(e^{\beta u^{2}}-1),\ \ \ \mbox{for} \ \mbox{all} \ u\in H^{1/2}(\mathbb{R}).
    \end{aligned}
\end{eqnarray*}
By Proposition $\ref{PRa}$, one has
\begin{eqnarray*}
    \begin{aligned}\displaystyle
    \int_{\mathbb{R}}(e^{\alpha|u|^{2}}-1)dx=\int_{\mathbb{R}}(e^{\alpha\|u\|^{2}(u/\|u\|)^{2}}-1)dx< +\infty, \ \mbox{for} \ \mbox{all} \ \|u\|\leq\sqrt{\pi/\alpha},
    \end{aligned}
\end{eqnarray*}
and
\begin{eqnarray*}
    \begin{aligned}\displaystyle
    \int_{\mathbb{R}}(e^{\beta|u|^{2}}-1)dx=\int_{\mathbb{R}}(e^{\beta\|u\|^{2}(u/\|u\|)^{2}}-1)dx< +\infty \ \mbox{for} \ \mbox{all} \ \|u\|\leq\sqrt{\pi/\beta}.
    \end{aligned}
\end{eqnarray*}
Then we deduce that
\begin{eqnarray*}
    \begin{aligned}\displaystyle
    \int_{\mathbb{R}}\Big(I_{\mu_{2}}\ast F(u)\Big)F(u)dx\leq|F(u)|_{\frac{2}{1+\mu_{2}}}|F(u)|_{\frac{2}{1+\mu_{2}}}\leq (\frac{\varepsilon}{2}\|u\|_{1/2}^{2}+C_{1}|u|_{3}^{3})^{2},
    \end{aligned}
    \end{eqnarray*}
and
\begin{eqnarray*}
    \begin{aligned}\displaystyle
    \int_{\mathbb{R}}\Big(I_{\mu_{1}}\ast G(u)\Big)G(u)dx&\leq|G(u)|_{\frac{2}{1+\mu_{1}}}|G(u)|_{\frac{2}{1+\mu_{1}}}\leq (\frac{\varepsilon}{2}\|u\|_{1/2}^{2}+C_{2}|u|_{3}^{3})^{2}.
    \end{aligned}
    \end{eqnarray*}
Together with the continuous embedding $H^{1/2}(\mathbb{R})\hookrightarrow L^{3}(\mathbb{R})$, taking $\sigma>0$, we obtain
 \begin{eqnarray*}
    \begin{aligned}\displaystyle
      \Phi(z)\geq (\frac{1}{2}-\frac{C\varepsilon}{2})\|u\|_{1/2}^{2}-C\|u\|_{1/2}^{6}\geq\sigma>0.
    \end{aligned}
 \end{eqnarray*}
 If we take $\rho$ to be sufficiently small, the the proof of the lemma is finished.
\end{proof}
Let $e\in H^{1/2}(\mathbb{R})\backslash \{0\}$ with $\|e\|_{1/2}=1$ and
\begin{eqnarray*}
    \begin{aligned}\displaystyle
Q=\{r(e,e)+(w,-w):w\in H^{1/2}(\mathbb{R}), \|w\|_{1/2}\leq R_{0} \ and \ 0\leq r\leq R_{1}\},
    \end{aligned}
\end{eqnarray*}
where $R_{0},R_{1}>\rho$ are defined in the Lemma $\ref{LEce}$.

\begin{lemma}\label{LEce}
Suppose that $(H_{0})-(H_{1})$ hold. Then there exist positive constants $R_{0}$, $R_{1}>\rho$ such that $\sup\Phi(\partial \mathcal{Q})\leq 0$.
\end{lemma}

\begin{proof}
Notice that the boundary $\partial Q$ of the set $Q$ is composed of three parts.

$(i)$ If $\omega\in \partial Q \cap E^{-}$, then $\omega=(u,u)\in E^{-}$, hence
\begin{eqnarray*}
    \begin{aligned}\displaystyle
\Phi(w)= -\|u\|^{2}_{1/2}-\frac{1}{2}\int_{\mathbb{R}}\Big(I_{\mu_{2}}\ast F(u)\Big)F(u)dx-\frac{1}{2}\int_{\mathbb{R}}\Big(I_{\mu_{1}}\ast G(-u)\Big)G(-u)dx \leq 0,
    \end{aligned}
\end{eqnarray*}
because $F,G$ are nonnegative functions.

$(ii)$ Suppose that $\omega=r(e,e)+(u,-u)\in\partial Q$ with $\|(u,-u)\|_{E}=R_{0}$ and $0\leq r\leq R_{1}$, we obtain
\begin{eqnarray*}
    \begin{aligned}\displaystyle
\Phi(\omega)=& r^{2}\|e\|^{2}_{1/2}-\|u\|^{2}_{1/2}-\frac{1}{2}\int_{\mathbb{R}}\Big(I_{\mu_{2}}\ast F(re+u)\Big)F(re+u)dx\\
&-\frac{1}{2}\int_{\mathbb{R}}\Big(I_{\mu_{1}}\ast G(re-u)\Big)G(re-u)dx
\leq R_{1}^{2}-\frac{1}{2}R_{0}^{2}.
    \end{aligned}
\end{eqnarray*}
Hence $\Phi(\omega)\leq0$ provided that $\sqrt{2}R_{1}\leq R_{0}$.

$(iii)$ If $\omega = R_{1}(e,e)+(u,-u)\in\partial Q$, with $\|(u,-u)\|_{E}\leq R_{0}$ for $R_{0}$ given by case $(ii)$. In this case,
\begin{eqnarray*}
    \begin{aligned}\displaystyle
    \Phi (\omega) =& R_{1}\|e\|^{2}_{1/2}-\|u\|^{2}_{1/2}-\frac{1}{2}\int_{\mathbb{R}}\Big(I_{\mu_{2}}\ast F(R_{1}e+u)\Big)F(R_{1}e+u)dx\\
    &-\frac{1}{2}\int_{\mathbb{R}}\Big(I_{\mu_{1}}\ast G(R_{1}e-u)\Big)G(R_{1}e-u)dx.
    \end{aligned}
\end{eqnarray*}

Define
\begin{eqnarray*}
    \begin{aligned}\displaystyle
    H(z)=\int_{\mathbb{R}}\Big(I_{\mu_{2}}\ast F(z)\Big)F(z)dx,
    \end{aligned}
\end{eqnarray*}
and
\begin{eqnarray*}
    \begin{aligned}\displaystyle
    s(t)= H\Big(t\frac{u/R_{1}+e}{\|u\|_{1/2}+\|e\|_{1/2}}\Big).
    \end{aligned}
\end{eqnarray*}

It follows from $(H_{2})$ that
\begin{eqnarray*}
    \begin{aligned}\displaystyle
    \frac{s^{\prime}(t)}{s(t)}\geq\frac{2\theta}{t},
    \end{aligned}
\end{eqnarray*}
where $\theta$ is defined in $(H_{2})$, which implies
\begin{eqnarray*}
    \begin{aligned}\displaystyle
    H(u+R_{1}e)=s(R_{1}(\|u\|_{1/2}+\|e\|_{1/2}))\geq C[R_{1}(\|u\|_{1/2}+\|e\|_{1/2})]^{2\theta}.
    \end{aligned}
\end{eqnarray*}
Then, we obtain
\begin{eqnarray*}
    \begin{aligned}\displaystyle
    \Phi (\omega)\leq R_{1}^{2}\|e\|^{2}_{1/2}-\frac{1}{2}H(u+R_{1}e)\leq R_{1}^{2}\|e\|^{2}_{1/2}-C R_{1}^{2\theta}\|y\|^{2\theta}_{1/2}\leq R_{1}^{2}-C R_{1}^{2\theta}.
    \end{aligned}
\end{eqnarray*}
Since $\theta>1$, taking $R_{1}$ sufficiently large, we obtain $\Phi(\omega)\leq0$.
\end{proof}

\section{{\bfseries Finite dimensional approximation}}\label{FIn}

In this section, since the functional $\Phi$ is strongly indefinite and defined in an infinite dimensional space, we use a finite dimensional approximation to produce a sequence approximating system $(\ref{a})$ of finite dimensional problems.

For the eigenvalues $0<\lambda_{1}<\lambda_{2}<\lambda_{3}<\cdot\cdot\cdot<\lambda_{j} \rightarrow +\infty$ of $(-\Delta +1,H^{1/2}(\mathbb{R}))$, there exists an orthonormal basis $\{\phi_{1},\phi_{2},...\}$ of corresponding eigenfunctions in $H^{1/2}(\mathbb{R})$. We set
\begin{eqnarray*}
    \begin{aligned}\displaystyle
    &E_{n}^{+}= span \{(\phi_{i},\phi_{i})|i=1,...,n\},\\
    &E_{n}^{-}= span \{(\phi_{i},\phi_{-i})|i=1,...,n\},\\
    &E_{n}=E_{n}^{+}\oplus E_{n}^{-}.
    \end{aligned}
\end{eqnarray*}
Define the hole functions (see \cite{Ca,Milan}) $\zeta_{m}: \mathbb{R}\rightarrow \mathbb{R}$ as
\begin{eqnarray*}
    \begin{aligned}\displaystyle
\zeta_{m}(x):=
    \left\{ \arraycolsep=1.5pt
       \begin{array}{ll}
        0 \ \ \ & 0\leq |x|\leq \frac{1}{m},\\[2mm]
        2+2 \frac{\log|x|}{\log m} \ \ \ & \frac{1}{m} <|x|< \frac{1}{\sqrt{m}},\\[2mm]
        1,\ \ \ & |x|\geq \frac{1}{\sqrt{m}}.
        \end{array}
    \right.
\end{aligned}
\end{eqnarray*}
For each $n,m \in \mathbb{N}$, consider the following finite-dimensional subspaces:
\begin{eqnarray*}
    \begin{aligned}\displaystyle
    &E_{n,m}:=\{u_{m}:=\zeta_{m}u : u\in E_{n}\},\\
    &E_{n,m}^{+}:=\{(u,u):u\in E_{n,m}\},\\
    &E_{n,m}^{-}:=\{(u,-u):u\in E_{n,m}\}.
    \end{aligned}
\end{eqnarray*}
Let $y\in H^{1/2}(\mathbb{R})$ be a fixed nonnegative function and
\begin{eqnarray*}
    \begin{aligned}\displaystyle
Q_{n,m,y}=\{r(y,y)+w:w\in E_{n,m}^{-}, \|w\|_{1/2}\leq R_{0} \ and \ 0\leq r\leq R_{1}\},
\end{aligned}
\end{eqnarray*}
where $R_{0},R_{1}$ are given in Lemma $\ref{LEce}$. We recall that these constants only depend of $y$. Then, we use the following notation:
\begin{eqnarray*}
    \begin{aligned}\displaystyle
H_{n,m,y}=\mathbb{R}(y,y)\oplus E_{n,m}, \ \ H_{n,m,y}^{+}=\mathbb{R}(y,y)\oplus E_{n,m}^{+}, \ \ H_{n,m,y}^{-}=\mathbb{R}(y,y)\oplus E_{n,m}^{-}.
    \end{aligned}
\end{eqnarray*}
Furthermore, define the class of mappings
\begin{eqnarray*}
    \begin{aligned}\displaystyle
\Gamma_{n,m,y}=\{h\in C(Q_{n,m,y},H_{n,m,y}): h(z)=z \ on \ \partial Q_{n,m,y}\}
    \end{aligned}
\end{eqnarray*}
and set
\begin{eqnarray*}
    \begin{aligned}\displaystyle
    c_{n,m,y}=\inf\limits_{h\in \Gamma_{n,m,y}} \max\limits_{z\in Q_{n,m,y}} \Phi(h(z)).
    \end{aligned}
\end{eqnarray*}

Using an intersection theorem (see Proposition $5.9$ in \cite{Rabi} ), we obtain
\begin{eqnarray*}
    \begin{aligned}\displaystyle
    h(Q_{n,m,y})\cap (\partial B_{\rho}\cap E^{+})\neq \emptyset, \ \forall h\in \Gamma_{n,m,y},
    \end{aligned}
\end{eqnarray*}
which in combination with Lemma $\ref{LEcd}$ implies that $c_{n,m,y}\geq \sigma>0$. On the other hand, since the identity mapping $Id:Q_{n,m,y}\rightarrow H_{n,m,y}$ belongs to $\Gamma_{n,m,y}$, we have for $z=r(y,y)+(u,-u)\in Q_{n,m,y}$ that
\begin{eqnarray*}
    \begin{aligned}\displaystyle
    \Phi(z)=&r^{2}\|y\|^{2}_{1/2}-\|u\|^{2}_{1/2}\\
    &-\frac{1}{2}\int_{\mathbb{R}}\Big(I_{\mu_{2}}\ast F(ry+u)\Big)F(ry+u)dx-\frac{1}{2}\int_{\mathbb{R}}\Big(I_{\mu_{1}}\ast G(ry-u)\Big)G(ry-u)dx\\
    \leq& r^{2}\|y\|_{1/2}^{2}\leq R_{1}^{2}.
    \end{aligned}
\end{eqnarray*}
Then we have
\begin{eqnarray*}
    \begin{aligned}\displaystyle
    0< \sigma \leq c_{n,m,y} \leq R_{1}^{2}.
    \end{aligned}
\end{eqnarray*}
We remark that the upper bound does not depend of $n$, but it depends on $y$.

Now, we denote $\Phi_{n,m,y}$ the restriction of the functional $\Phi$ to the finite dimensional subspace $H_{n,m,y}$. Therefore, in view of Lemmas $\ref{LEcd}$ and $\ref{LEce}$, we see that the geometry of a linking theorem holds for the functional $\Phi_{n,m,y}$. Thus, applying the linking theorem for $\Phi_{n,m,y}$ (see Theorem $5.3$ in \cite{Rabi}), we obtain a $(PS)$ sequence, which is bounded in view of Lemma $\ref{LEcb}$. Finally, using the fact that $H_{n,m,y}$ is a finite dimensional space, we get the main result of this section.
\begin{proposition}\label{PRda}
For each $n,m\in \mathbb{N}$ and for each $y\in H^{1/2}(\mathbb{R})$, a fixed nonnegative function, the functional $\Phi_{n,m,y}$ has a critical point at level $c_{n,m,y}$. More precisely, there is a $z_{n,m,y}\in H_{n,m,y}$ such that
\begin{eqnarray*}
    \begin{aligned}\displaystyle
    \Phi_{n,m,y}(z_{n,m,y})=c_{n,m,y}\in [\sigma, R_{1}^{2}] \ \ and \ \ (\Phi_{n,m,y})^{\prime}(z_{n,m,y})=0.
    \end{aligned}
\end{eqnarray*}
Furthermore, $\|z_{n,m,y}\|_{1/2}\leq C$ where $C$ does not depend of $n$.
\end{proposition}

\section{{\bfseries The estimates for the critical level}}\label{TECL}
In this section, we assume that $f$ and $g$ have critical growth with exponent critical $\alpha_{0}$ and $\beta_{0}$. Then together Trudinger-Moser inequality with $(H_{5})$, we obtain an upper bound for the minimax level. In order to do this, we prove the following result in which we combine the hole functions, the Morser type functions and an approximation argument inspired by Tang et al. \cite{QTZA}.

Let us introduce the following Moser type functions (see \cite{Do}) supported in $B_{1}$ given by
\begin{eqnarray*}
    \begin{aligned}\displaystyle
    w_{m}(x)=\frac{1}{\sqrt{\pi}}
    \left\{ \arraycolsep=1.5pt
       \begin{array}{ll}
        \sqrt{\log m} \ \ \ & 0\leq |x|\leq 1/m;\\[2mm]
        \frac{\log \frac{1}{|x|}}{\sqrt{\log m}} \ \ \ & 1/m \leq|x|\leq 1;\\[2mm]
        0,\ \ \ & |x|\geq 1.
        \end{array}
    \right.
    \end{aligned}
\end{eqnarray*}
One has that
\begin{eqnarray*}
    \begin{aligned}\displaystyle
    \|(-\Delta)^{\frac{1}{4}} w_{m}\|^{2}_{2}=\int_{\mathbb{R}}|(-\Delta)^{\frac{1}{4}} w_{m}|^{2}dx =1+o(1), \ \ \ \mbox{and} \ \ \| w_{m}\|^{2}_{2}=\int_{\mathbb{R}}| w_{m}|^{2}dx=\delta_{m},
    \end{aligned}
\end{eqnarray*}
where
\begin{equation}\label{ESTa}
    \begin{aligned}\displaystyle
    \delta_{m} &=\frac{1}{\pi}\Big(\int^{\frac{1}{m}}_{-\frac{1}{m}}\log m \ dx +\frac{1}{\log m}\int^{-\frac{1}{m}}_{-1}(\log |x|)^{2}dx+\frac{1}{\log m}\int^{1}_{\frac{1}{m}}(\log |x|)^{2}dx\Big)\\
    &=\frac{4}{\pi}\Big(\frac{1}{\log m}-\frac{1}{m\log m}-\frac{1}{m}\Big)> 0.
    \end{aligned}
\end{equation}
Thanks to $(H_{5})$ we have for all $t \geq R_{\varepsilon}$,\\
\begin{equation}\label{ESTb}
    \begin{aligned}\displaystyle
tF(t)\geq (\kappa_{1}-\varepsilon) e^{\alpha_{0}t^{2}}\ \ and \ \ tG(t)\geq (\kappa_{2}-\varepsilon) e^{\beta_{0}t^{2}}.
    \end{aligned}
\end{equation}

\begin {lemma}\label{LEfa}Suppose that $(H_{0}),(H_{1}),(H_{4})$ and $(H_{5})$ hold, then there exists $m_{0}\in  \mathbb{N}$ such that for all $m\geq m_{0}$ the corresponding Moser's function $w_{m}$ satisfies
\begin{equation}\label{ESTc}
    \begin{aligned}\displaystyle
    \sup\limits_{\mathbb{R^{+}}(w_{m},w_{m})\bigoplus \mathbb{E^{-}}} \Phi< \min\{ \frac{(1+\mu_{1})\pi}{2\beta_{0}},\frac{(1+\mu_{2})\pi}{2\alpha_{0}}\}.
    \end{aligned}
\end{equation}
\end{lemma}

\begin{proof}
First, we give the following estimate in $B_{1/m}(0)$,
\begin{eqnarray*}
    \begin{aligned}\displaystyle
\int^{\frac{1}{m}}_{-\frac{1}{m}}\int^{\frac{1}{m}}_{-\frac{1}{m}}\frac{dxdy}{|x-y|^{1-\mu}}=\frac{2^{2+\mu}}{\mu(1+\mu)}\Big(\frac{1}{m}\Big)^{1+\mu}:=C_{\mu}\Big(\frac{1}{m}\Big)^{1+\mu}.
    \end{aligned}
\end{eqnarray*}
By $(H_{5})$, we may choose $\varepsilon>0$ small such that
\begin{equation}\label{import}
    \begin{aligned}\displaystyle
    \log \frac{2A_{\mu_{1}}C_{\mu_{1}}(\kappa_{2}-\varepsilon)^{2}\beta_{0}^{2}}{\pi  (1+\varepsilon)^{2}(1+\mu_{1})}>\frac{4(1+\mu_{1})}{\pi}-1
    \end{aligned}
\end{equation}
and
\begin{equation}
    \begin{aligned}\displaystyle
    \log \frac{2A_{\mu_{2}}C_{\mu_{2}}(\kappa_{1}-\varepsilon)^{2}\alpha_{0}^{2}}{\pi  (1+\varepsilon)^{2}(1+\mu_{2})}>\frac{4(1+\mu_{2})}{\pi}-1.
    \end{aligned}
\end{equation}
From Proposition $\ref{PRda}$ and $(\ref{ESTa})$, we have
\begin{eqnarray*}
    \begin{aligned}\displaystyle
    \Phi(tw_{m}+v_{m},tw_{n}-v_{m})
    &=t^{2}\|w_{m}\|^{2}_{1/2}-\|v_{m}\|^{2}_{1/2}-\Psi(tw_{m}+v_{m},tw_{m}-v_{m})\\
    &\leq t^{2}(\|(-\Delta)^{\frac{1}{4}} w_{m}\|^{2}_{2}+\|w_{m}\|^{2}_{2})-\Psi(tw_{m}+v_{m},tw_{m}-v_{m})\\
    &\leq (1+\delta_{m})t^{2}-\Psi(tw_{m}+v_{m},tw_{m}-v_{m}),\ \ \ \forall \ t\geq 0.
    \end{aligned}
\end{eqnarray*}
where $\Psi(tw_{m}+v_{m},tw_{m}-v_{m})=\frac{1}{2}\int_{\mathbb{R}}(I_{\mu_{2}}\ast F(tw_{m}+v_{m}))F(tw_{m}+v_{m})dx+\frac{1}{2}\int_{\mathbb{R}}(I_{\mu_{1}}\ast G(tw_{m}-v_{m}))G(tw_{m}-v_{m})dx$, and $u_{m}\in E_{n,m}$.\\
When $\frac{(1+\mu_{1})\pi}{2\beta_{0}}\leq\frac{(1+\mu_{2})\pi}{2\alpha_{0}}$, there are four cases to distinguish.

Case $i).$ $t\in[0,\sqrt{\frac{(1+\mu_{1})\pi}{4\beta_{0}}}]$. Then by $F,G\geq0$, we have
\begin{eqnarray*}
    \begin{aligned}\displaystyle
    \Phi(tw_{n}+v,tw_{n}-v)\leq(1+\frac{4}{\pi\log n})t^{2}-\Psi(tw_{n}+v,tw_{n}-v)\leq \frac{(1+\mu_{1})\pi}{4\beta_{0}}+O(\frac{1}{\log m}),
    \end{aligned}
\end{eqnarray*}
which implies that $(\ref{ESTc})$ hold.

Case $ii).$ $t\in[\sqrt{\frac{(1+\mu_{1})\pi}{4\beta_{0}}},\sqrt{\frac{(1+\mu_{1})\pi}{2\beta_{0}}}]$. In this case, $tw_{m}\geq R_{\varepsilon}$ for $x\in B_{1/m}(0)$ and $m\in \mathbb{N}$ large. Then using $(3.1)$, we have\\
\begin{equation}\label{Newed}
    \begin{aligned}\displaystyle
    F(tw_{m}+v_{m})\geq\frac{(\kappa_{1}-\varepsilon)e^{\alpha_{0}(tw_{m}+v_{m})^{2}}}{tw_{m}+v_{m}} \ \ \ and \ \ \ \ G(tw_{m}-v_{m})\geq\frac{(\kappa_{2}-\varepsilon)e^{\beta_{0}(tw_{m}-v_{m})^{2}}}{tw_{m}-v_{m}},
    \end{aligned}
\end{equation}
which, together with $(\ref{ESTb})$, yields
\begin{eqnarray*}
    \begin{aligned}\displaystyle
    &\Psi(tw_{m}+v_{m},tw_{m}-v_{m})\\
    &=\frac{1}{2}\int_{\mathbb{R}}\Big(I_{\mu_{2}}\ast F(tw_{m}+v_{m})\Big)F(tw_{m}+v_{m})dx+\frac{1}{2}\int_{\mathbb{R}}\Big(I_{\mu_{1}}\ast G(tw_{m}-v_{m})\Big)G(tw_{m}-v_{m})dx\\
    &\geq \frac{A_{\mu_{2}}}{2}\int_{-\frac{1}{m}}^{\frac{1}{m}}\int_{-\frac{1}{m}}^{\frac{1}{m}}\frac{F(tw_{m}+v_{m}(x))\cdot F(tw_{m}+v_{m}(y))}{|x-y|^{1-\mu_{2}}}dxdy\\
    &\quad+\frac{A_{\mu_{1}}}{2}\int_{-\frac{1}{m}}^{\frac{1}{m}}\int_{-\frac{1}{m}}^{\frac{1}{m}}\frac{G(tw_{m}-v_{m}(x))\cdot G(tw_{m}-v_{m}(y))}{|x-y|^{1-\mu_{1}}}dxdy.
    \end{aligned}
\end{eqnarray*}
Using the functions $F,G$ are nonnegative and $u_{m}$ is zero for all $|x|\leq\frac{1}{m}$, we have
\begin{equation}\label{Newee}
    \begin{aligned}\displaystyle
    \Psi(tw_{m}+v_{m},tw_{m}-v_{m})\geq \frac{A_{\mu_{2}}}{2}\int_{-\frac{1}{m}}^{\frac{1}{m}}\int_{-\frac{1}{m}}^{\frac{1}{m}}\frac{F(tw_{m})\cdot F(tw_{m})}{|x-y|^{1-\mu_{2}}}dxdy,
    \end{aligned}
\end{equation}
and
\begin{equation}\label{Newef}
    \begin{aligned}\displaystyle
    \Psi(tw_{m}+v_{m},tw_{m}-v_{m})\geq
    \frac{A_{\mu_{1}}}{2}\int_{-\frac{1}{m}}^{\frac{1}{m}}\int_{-\frac{1}{m}}^{\frac{1}{m}}\frac{G(tw_{m})\cdot G(tw_{m})}{|x-y|^{1-\mu_{1}}}dxdy.
    \end{aligned}
\end{equation}
Next by $(\ref{Newed})$ and $(\ref{Newef})$, we obtain
\begin{eqnarray*}
    \begin{aligned}\displaystyle
    \Psi(tw_{m}+v_{m},tw_{m}-v_{m})
    &\geq \frac{A_{\mu_{1}}(\kappa_{2}-\varepsilon)^{2}}{2t^{2}w_{m}^{2}}\int_{-\frac{1}{m}}^{\frac{1}{m}}\int_{-\frac{1}{m}}^{\frac{1}{m}}\frac{e^{2\beta_{0}t^{2}w_{m}^{2}}}{|x-y|^{1-\mu_{1}}}dxdy\\
    &\geq\frac{A_{\mu_{1}}(\kappa_{2}-\varepsilon)^{2}\beta_{0}}{(1+\mu_{1})\log m}e^{2\beta_{0}t^{2}w_{m}^{2}}\int_{-\frac{1}{m}}^{\frac{1}{m}}\int_{-\frac{1}{m}}^{\frac{1}{m}}\frac{1}{|x-y|^{1-\mu_{1}}}dxdy\\
    &\geq \frac{A_{\mu_{1}}(\kappa_{2}-\varepsilon)^{2}\beta_{0}C_{\mu_{1}}}{(1+\mu_{1})\log m}\Big(\frac{1}{m}\Big)^{1+\mu_{1}}e^{2\beta_{0}t^{2}w_{m}^{2}}\\
    &\geq \frac{A_{\mu_{1}}C_{\mu_{1}}(\kappa_{2}-\varepsilon)^{2}\beta_{0}}{(1+\mu_{1}) m^{1+\mu_{1}}\log m}e^{2(\pi)^{-1}\beta_{0}t^{2}\log m}.
    \end{aligned}
\end{eqnarray*}
Then we obtain that
\begin{eqnarray*}
    \begin{aligned}\displaystyle
    &\Phi(tw_{m}+v_{m},tw_{m}-v_{m})\leq (1+\frac{4}{\pi\log m})t^{2} -\frac{A_{\mu_{1}}C_{\mu_{1}}(\kappa_{2}-\varepsilon)^{2}\beta_{0}}{(1+\mu_{1}) m^{1+\mu_{1}}\log m}e^{2(\pi)^{-1}\beta_{0}t^{2}\log m}:=\varphi_{m}(t).
    \end{aligned}
\end{eqnarray*}
Choosing $t_{m}>0$ satisfying $\varphi_{m}^{\prime}(t_{m})=0$, then we obtain
\begin{eqnarray*}
    \begin{aligned}\displaystyle
    1+\frac{4}{\pi\log m}=\frac{2A_{\mu_{1}}C_{\mu_{1}}(\kappa_{2}-\varepsilon)^{2}\beta_{0}^{2}}{ (1+\mu_{1}) \pi m^{1+\mu_{1}}}e^{2(\pi)^{-1}\beta_{0}t^{2}\log m}.
    \end{aligned}
\end{eqnarray*}
It follows that
\begin{equation}\label{ESTd}
    \begin{aligned}\displaystyle
    t_{m}^{2}&=\frac{(1+\mu_{1})\pi}{2\beta_{0}}\Big[1+\frac{\log(1+\frac{4}{\pi\log m})-\log\frac{2A_{\mu_{1}}C_{\mu_{1}}(\kappa_{2}-\varepsilon)^{2}\beta_{0}^{2}}{\pi (1+\mu_{1})}}{(1+\mu_{1})\log m}\Big]\\
    &\leq \frac{(1+\mu_{1})\pi}{2\beta_{0}}-\frac{\pi}{2\beta_{0}\log m}\log \frac{2A_{\mu_{1}}C_{\mu_{1}}(\kappa_{2}-\varepsilon)^{2}\beta_{0}^{2}}{\pi (1+\mu_{1})  (1+\varepsilon )},
    \end{aligned}
\end{equation}
and
\begin{eqnarray*}
    \begin{aligned}\displaystyle
    \varphi_{m}(t)\leq \varphi_{m}(t_{m})=(1+\frac{4}{\pi\log m})t_{m}^{2} -\frac{\pi }{2\beta_{0}\log m}(1+\frac{4}{\pi\log m}), \ \ \forall \ t\geq 0.
    \end{aligned}
\end{eqnarray*}
Using $(\ref{ESTd})$, we get
\begin{eqnarray*}
    \begin{aligned}\displaystyle
    \varphi_{m}(t)&\leq (1+\frac{4}{\pi\log m})t_{m}^{2}-\frac{\pi }{2\beta_{0}\log m}(1+\frac{4}{\pi\log m})\\
    \leq& (1+\frac{4}{\pi\log m})\Big[\frac{(1+\mu_{1})\pi}{2\beta_{0}}-\frac{\pi}{2\beta_{0}\log m}\log \frac{2A_{\mu_{1}}C_{\mu_{1}}(\kappa_{2}-\varepsilon)^{2}\beta_{0}^{2}}{\pi (1+\mu_{1})  (1+\varepsilon )}\Big]-\frac{\pi}{2\beta_{0}\log m}\\&+O(\frac{1}{\log^{2} m})\\
    \leq& \frac{(1+\mu_{1})\pi}{2\beta_{0}}+\frac{2(1+\mu_{1})}{\beta_{0}\log m}-\frac{\pi}{2\beta_{0}\log m}-\frac{\pi}{2\beta_{0}\log m}\log \frac{2A_{\mu_{1}}C_{\mu_{1}}(\kappa_{2}-\varepsilon)^{2}\beta_{0}^{2}}{\pi (1+\mu_{1})  (1+\varepsilon )}+O(\frac{1}{\log^{2} m}).
    \end{aligned}
\end{eqnarray*}
Hence, we obtain
\begin{eqnarray*}
    \begin{aligned}\displaystyle
    \Phi(tw_{m}+v_{m},tw_{m}-v_{m})\leq& \frac{(1+\mu_{1})\pi}{2\beta_{0}}+\frac{\pi}{2\beta_{0}\log m}\Big[\frac{4(1+\mu_{1})}{\pi}-1-\log \frac{2A_{\mu_{1}}C_{\mu_{1}}(\kappa_{2}-\varepsilon)^{2}\beta_{0}^{2}}{\pi (1+\mu_{1})  (1+\varepsilon )}\Big]\\
    &+O(\frac{1}{\log^{2} m}),
    \end{aligned}
\end{eqnarray*}
which together with $(\ref{import})$ imply that $(\ref{ESTc})$ hold.

Case $iii).$ $t\in\Big[\sqrt{\frac{(1+\mu_{1})\pi}{2\beta_{0}}},\sqrt{\frac{(1+\mu_{1})\pi}{2\beta_{0}}(1+\varepsilon)}\Big]$. In this case, $tw_{m}\geq R_{\varepsilon}$ for $x\in B_{1/m}(0)$ and $m\in \mathbb{N}$ large. Then using $(\ref{ESTb})$, we have
\begin{eqnarray*}
    \begin{aligned}\displaystyle
    \Psi(tw_{m}+v_{m},tw_{m}-v_{m})&\geq \frac{A_{\mu_{1}}(\kappa_{2}-\varepsilon)^{2}}{2t^{2}w_{m}^{2}}\int_{-\frac{1}{m}}^{\frac{1}{m}}\int_{-\frac{1}{m}}^{\frac{1}{m}}\frac{e^{2\beta_{0}t^{2}w_{m}^{2}}}{|x-y|^{1-\mu_{1}}}dxdy\\
    &\geq \frac{A_{\mu_{1}}C_{\mu_{1}}(\kappa_{2}-\varepsilon)^{2}\beta_{0}}{(1+\varepsilon)(1+\mu_{1}) m^{1+\mu_{1}}\log m}e^{2(\pi)^{-1}\beta_{0}t^{2}\log m}.
    \end{aligned}
\end{eqnarray*}
Then
\begin{eqnarray*}
    \begin{aligned}\displaystyle
    \Phi(tw_{m}+v_{m},tw_{m}-v_{m})&\leq (1+\frac{4}{\pi\log m})t^{2} -\frac{A_{\mu_{1}}C_{\mu_{1}}(\kappa_{2}-\varepsilon)^{2}\beta_{0}}{(1+\varepsilon)(1+\mu_{1}) m^{1+\mu_{1}}\log m}e^{2(\pi)^{-1}\beta_{0}t^{2}\log m}\\
    &:=\psi_{m}(t).
    \end{aligned}
\end{eqnarray*}
Choosing $\widetilde{t}_{m}>0$ satisfying $\psi_{m}^{\prime}(\widetilde{t}_{m})=0$, then we obtain
\begin{eqnarray*}
    \begin{aligned}\displaystyle
    1+\frac{4}{\pi\log m}=\frac{2A_{\mu_{1}}C_{\mu_{1}}(\kappa_{2}-\varepsilon)^{2}\beta_{0}^{2}}{ (1+\varepsilon)(1+\mu_{1}) \pi m^{1+\mu_{1}}}e^{2(\pi)^{-1}\beta_{0}t^{2}\log m}.
    \end{aligned}
\end{eqnarray*}
It follows that
\begin{equation}\label{ESTe}
    \begin{aligned}\displaystyle
    \widetilde{t}_{m}^{2}&=\frac{(1+\mu_{1})\pi}{2\beta_{0}}\Big[1+\frac{\log(1+\frac{4}{\pi\log m})-\log\frac{2A_{\mu_{1}}C_{\mu_{1}}(\kappa_{2}-\varepsilon)^{2}\beta_{0}^{2}}{(1+\varepsilon)\pi (1+\mu_{1})}}{(1+\mu_{1})\log m}\Big]\\
    &\leq \frac{(1+\mu_{1})\pi}{2\beta_{0}}+\frac{\pi}{2\beta_{0}\log m}\log \frac{(1+\varepsilon)\pi (1+\mu_{1})  (1+\varepsilon )}{2A_{\mu_{1}}C_{\mu_{1}}(\kappa_{2}-\varepsilon)^{2}\beta_{0}^{2}},
    \end{aligned}
\end{equation}
and
\begin{eqnarray*}
    \begin{aligned}\displaystyle
    \psi_{m}(t)\leq \psi_{m}(\widetilde{t}_{m})=(1+\frac{4}{\pi\log m})\widetilde{t}_{m}^{2} -\frac{\pi }{2\beta_{0}\log m}(1+\frac{4}{\pi\log m}), \ \ \forall \ t \ \geq \ 0.
    \end{aligned}
\end{eqnarray*}
Using $(\ref{ESTe})$, we get
\begin{eqnarray*}
    \begin{aligned}\displaystyle
    \psi_{m}(t)&\leq (1+\frac{4}{\pi\log m})\widetilde{t}_{m}^{2} -\frac{\pi }{2\beta_{0}\log m}(1+\frac{4}{\pi\log m})\\
    \leq& (1+\frac{4}{\pi\log m})\Big[\frac{(1+\mu_{1})\pi}{2\beta_{0}}+\frac{\pi}{2\beta_{0}\log m}\log \frac{(1+\varepsilon)^{2}\pi (1+\mu)}{2A_{\mu_{1}}C_{\mu_{1}}(\kappa_{2}-\varepsilon)^{2}\beta_{0}^{2}}\Big]-\frac{\pi }{2\beta_{0}\log m}\\&+O(\frac{1}{\log^{2} m})\\
    \leq& \frac{(1+\mu_{1})\pi}{2\beta_{0}} +\frac{\pi}{2\beta_{0}\log m}\Big[ \frac{4(1+\mu_{1})}{\pi}-1-\log \frac{2A_{\mu_{1}}C_{\mu_{1}}(\kappa_{2}-\varepsilon)^{2}\beta_{0}^{2}}{(1+\varepsilon)^{2}\pi (1+\mu_{1})  }\Big]+O(\frac{1}{\log^{2} m}).
    \end{aligned}
\end{eqnarray*}
Hence, we obtain
\begin{eqnarray*}
    \begin{aligned}\displaystyle
    \Phi(tw_{m}+v_{m},tw_{m}-v_{m})\leq& \frac{(1+\mu_{1})\pi}{2\beta_{0}}+\frac{\pi}{2\beta_{0}\log m}\Big[ \frac{4(1+\mu_{1})}{\pi}-1-\log \frac{2A_{\mu_{1}}C_{\mu_{1}}(\kappa_{2}-\varepsilon)^{2}\beta_{0}^{2}}{(1+\varepsilon)^{2}\pi (1+\mu_{1})} \Big]\\
    &+O(\frac{1}{\log^{2} m}),
    \end{aligned}
\end{eqnarray*}
which together with $(\ref{import})$ imply that $(\ref{ESTc})$ hold.

Case $iv).$ $t\in\Big(\sqrt{\frac{(1+\mu_{1})\pi}{2\beta_{0}}(1+\varepsilon)},+\infty\Big)$. In this case, $tw_{m}\geq R_{\varepsilon}$ for $x\in B_{1/m}(0)$ and $m\in \mathbb{N}$ large. Then using $(\ref{ESTb})$, we have
\begin{eqnarray*}
    \begin{aligned}\displaystyle
    \Phi(tw_{m}+v_{m},tw_{m}-v_{m})
    &\leq(1+\frac{4}{\pi\log m})t^{2}-\Psi(tw_{m}+v_{m},tw_{m}-v_{m})\\
    &\leq(1+\frac{4}{\pi\log m})t^{2}-\frac{\pi A_{\mu_{1}} C_{\mu_{1}} (\kappa_{2}-\varepsilon)^{2}}{  2t^{2} m^{1+\mu}\log m}e^{2(\pi)^{-1}\beta_{0}t^{2}\log m}\\
    &\leq (1+\frac{4}{\pi\log m})\frac{(1+\varepsilon)(1+\mu_{1})\pi}{2\beta_{0}}-\frac{\pi A_{\mu_{1}} C_{\mu_{1}} (\kappa_{2}-\varepsilon)^{2}}{  2t^{2}\log m}e^{(1+\mu_{1})\varepsilon\log m}\\
    &\leq (1+\frac{4}{\pi\log m})\frac{(1+\varepsilon)(1+\mu_{1})\pi}{2\beta_{0}}-\frac{ A_{\mu_{1}} C_{\mu_{1}} (\kappa_{2}-\varepsilon)^{2}\beta_{0}}{  (1+\mu_{1})\log m(1+\varepsilon)}e^{(1+\mu_{1})\varepsilon\log m}\\
    &< \frac{(1+\mu_{1})\pi}{2\beta_{0}},
    \end{aligned}
\end{eqnarray*}
where we have used the fact that the function
\begin{eqnarray*}
    \begin{aligned}\displaystyle
\phi(t)=(1+\frac{4}{\pi\log m})t^{2}-\frac{\pi A_{\mu_{1}} C_{\mu_{1}} (\kappa_{2}-\varepsilon)^{2}}{  2t^{2} m^{1+\mu}\log m}e^{2(\pi)^{-1}\beta_{0}t^{2}\log m}
    \end{aligned}
\end{eqnarray*}
is decreasing on $t\in\Big(\sqrt{\frac{(1+\mu_{1})\pi}{2\beta_{0}}(1+\varepsilon)},+\infty\Big)$. Thus $(\ref{ESTc})$ holds for $n\geq n_{0}$.\\
When $\frac{(1+\mu_{2})\pi}{2\alpha_{0}} \leq \frac{(1+\mu_{1})\pi}{2\beta_{0}}$, using the fact $u_{m}$ is zero for all $|x|\leq\frac{1}{m}$, we obtian
\begin{eqnarray*}
    \begin{aligned}\displaystyle
    \Psi(tw_{m}+v_{m},tw_{m}-v_{m})\geq \frac{A_{\mu_{2}}}{2}\int_{-\frac{1}{m}}^{\frac{1}{m}}\int_{-\frac{1}{m}}^{\frac{1}{m}}\frac{F(tw_{m})\cdot F(tw_{m})}{|x-y|^{1-\mu_{2}}}dxdy.
    \end{aligned}
\end{eqnarray*}
Following the same steps as the above, $(\ref{ESTc})$ is satisfied. Then the proof is complete.
\end{proof}

\section{{\bfseries Proof of Theorem 1.1.}}\label{PROOF}

In this section we assume that $f$ and $g$ have critical growth with exponent critical $\alpha_{0}$ and $\beta_{0}$.\\
{\bf Proof of Theorem} $\ref{THb}$. The Lemma $\ref{LEfa}$ implies that there is $\delta>0$ such that
\begin{eqnarray*}
    \begin{aligned}\displaystyle
    c_{n,m}:=c_{n,m,y}\leq \min\{ \frac{(1+\mu_{1})\pi}{2\beta_{0}},\frac{(1+\mu_{2})\pi}{2\alpha_{0}}\}-\delta:=Q\pi-\delta,
\end{aligned}
\end{eqnarray*}
where $Q=\min\{ \frac{(1+\mu_{1})}{2\beta_{0}},\frac{(1+\mu_{2})}{2\alpha_{0}}\}$. Applying Proposition $\ref{PRda}$ we have a sequence $z_{n,m}:=z_{n,m,y}=(u_{n,m},v_{n,m})\in H_{n,m,y}$ such that
\begin{equation}\label{PRCa}
    \begin{aligned}\displaystyle
    &\|(u_{n,m},v_{n,m})\|_{E}\leq C,\\
    &\Phi_{n,m,y}(u_{n,m},v_{n,m})=c_{n,m}:=c_{n,m,y}\in [\sigma,Q\pi-\delta),\\
    &(\Phi_{n,m,y})^{\prime}(u_{n,m},v_{n,m})=0,\\
    &(u_{n,m},v_{n,m})\rightharpoonup (u_{0},v_{0}) \ \ \mbox{in} \ \ E.
    \end{aligned}
\end{equation}

By $(u_{n,m},v_{n,m})$ is bounded in $E$, there exists $(u_{0},v_{0})$ such that $(u_{n,m},v_{n,m})\rightharpoonup (u_{0},v_{0})$.
Then Taking $(0,\psi)$ and $(\varphi,0)$ as text function in $(\ref{PRCa})$, we get
\begin{eqnarray*}
    \begin{aligned}\displaystyle
    &\int_{\mathbb{R}}(-\Delta)^{\frac{1}{4}} u_{n,m} (-\Delta)^{\frac{1}{4}} \psi +u_{n,m}\psi dx = \int_{\mathbb{R}}\Big(I_{\mu_{1}}\ast G(v_{n,m})\Big)g(v_{n,m})\psi dx \ \ \forall \psi \in E_{n,m}, \\
    &\int_{\mathbb{R}}(-\Delta)^{\frac{1}{4}} v_{n,m}(-\Delta)^{\frac{1}{4}} \varphi +v_{n,m}\varphi dx = \int_{\mathbb{R}}\Big(I_{\mu_{2}}\ast F(u_{n,m})\Big)f(u_{n,m})\varphi dx \ \ \forall \varphi \in E_{n,m},
    \end{aligned}
\end{eqnarray*}
where $E_{n,m}$ defined in Section $\ref{TECL}$. By Lemma $\ref{LEcb}$, we have
\begin{eqnarray*}
    \begin{aligned}\displaystyle
    &\lim\limits_{n\rightarrow\infty}\int_{\mathbb{R}}\Big(I_{\mu_{2}}\ast F(u_{n,m})\Big)f(u_{n,m})\varphi dx=\int_{\mathbb{R}}\Big(I_{\mu_{2}}\ast F(u_{0})\Big)f(u_{0})\varphi dx,\\
    &\lim\limits_{n\rightarrow\infty}\int_{\mathbb{R}}\Big(I_{\mu_{1}}\ast G(v_{n,m})\Big)g(v_{n,m})\psi dx=\int_{\mathbb{R}}\Big(I_{\mu_{1}}\ast G(v_{0})\Big)g(v_{0})\psi dx.
    \end{aligned}
\end{eqnarray*}
Thus, using the fact that $\cup_{n\in \mathbb{N}}E_{n,m}$ is dense in $E$, we can obtain that
\begin{equation}\label{PRCb}
    \begin{aligned}\displaystyle
    &\int_{\mathbb{R}}(-\Delta)^{\frac{1}{4}} u_{0}(-\Delta)^{\frac{1}{4}} \psi +u_{0}\psi dx = \int_{\mathbb{R}}\Big(I_{\mu_{1}}\ast G(v_{0})\Big)g(v_{0})\psi dx \ \ \forall \psi \in E, \\
    &\int_{\mathbb{R}}(-\Delta)^{\frac{1}{4}} v_{0}(-\Delta)^{\frac{1}{4}} \varphi +v_{0}\varphi dx = \int_{\mathbb{R}}\Big(I_{\mu_{2}}\ast F(u_{0})\Big)f(u_{0})\varphi dx \ \ \forall \varphi \in E.
    \end{aligned}
\end{equation}
Therefore, we conclude that $\Phi^{\prime}(u_{0},v_{0})=0$ in E, then $(u_{0},v_{0})$ is a critical point of $\Phi$. Then for $\forall \psi, \varphi \in C^{\infty}_{0}(\mathbb{R})$, we know that
\begin{equation}\label{PRCc}
    \begin{aligned}\displaystyle
    \langle u_{0}, \psi\rangle_{1/2} + \langle v_{0}, \varphi\rangle_{1/2} = &\int_{\mathbb{R}}\Big(I_{\mu_{2}}\ast F(u_{0})\Big)f(u_{0})\varphi dx+\int_{\mathbb{R}}\Big(I_{\mu_{1}}\ast G(v_{0})\Big)g(v_{0})\psi dx.
    \end{aligned}
\end{equation}

Now, it remains to prove that $u_{0}$, $v_{0}\neq 0$. Assume that $u_{0}\equiv 0$, then by $(\ref{PRCb})$ we know that $v_{0}=0$. Thus the proof of this theorem is divided into two steps.

Step $1.$ $(u_{n,m},v_{n,m})\rightarrow (u_{0},v_{0})=(0,0)$, that is, $\|u_{n,m}\|_{1/2},\|v_{n,m}\|_{1/2}\rightarrow 0$. By the Cauchy-Schwarz inequality, we obtain
\begin{eqnarray*}
    \begin{aligned}\displaystyle
    \lim\limits_{n\rightarrow \infty}\int_{\mathbb{R}}((-\Delta)^{\frac{1}{4}} u_{n,m}(-\Delta)^{\frac{1}{4}} v_{n,m}+ u_{n,m}v_{n,m})\rightarrow 0.
    \end{aligned}
\end{eqnarray*}
Which implies that
\begin{eqnarray*}
    \begin{aligned}\displaystyle
    \int_{\mathbb{R}}\Big(I_{\mu_{2}}\ast F(u_{n,m})\Big)f(u_{n,m})u_{n,m}dx\rightarrow 0 \ \ \mbox{and} \ \ \int_{\mathbb{R}}\Big(I_{\mu_{1}}\ast G(v_{n,m})\Big)g(v_{n,m})v_{n,m}dx\rightarrow 0,
    \end{aligned}
\end{eqnarray*}
together with $(H_{3})$ we get that
\begin{eqnarray*}
    \begin{aligned}\displaystyle
    \int_{\mathbb{R}}\Big(I_{\mu_{2}}\ast F(u_{n,m})\Big)F(u_{n,m})dx\rightarrow 0 \ \ \mbox{and} \ \  \int_{\mathbb{R}}\Big(I_{\mu_{1}}\ast G(v_{n,m})\Big)G(v_{n,m})dx\rightarrow 0.
    \end{aligned}
\end{eqnarray*}
Then, the last convergence shows us that $c_{n,m}=0$, a contradiction. Hence, this case cannot occur.

Step $2.$ $(u_{n,m},v_{n,m})$ converges weakly to $(u_{0},v_{0})$ in $E$ but does not converge strongly. In other words, $(u_{n,m},v_{n,m})\rightharpoonup (0,0)$ in $E$, and there is a constant $a>0$ such that $\liminf\limits_{n\rightarrow\infty}\|u_{n,m}\|_{1/2}\geq a$ and $\liminf\limits_{n\rightarrow\infty}\|v_{n,m}\|_{1/2}\geq a$.

Taking $(0,u_{n,m}),(v_{n,m},0)$ as test function in $(\ref{PRCa})$, we obtain that
\begin{eqnarray*}
    \begin{aligned}\displaystyle
    \|u_{n,m}\|_{1/2}^{2}=\int_{\mathbb{R}}\Big(I_{\mu_{1}}\ast G(v_{n,m})\Big)g(v_{n,m})u_{n,m}dx,
        \end{aligned}
\end{eqnarray*}
and
\begin{eqnarray*}
    \begin{aligned}\displaystyle
\ \|v_{n,m}\|_{1/2}^{2}=\int_{\mathbb{R}}\Big(I_{\mu_{2}}\ast F(u_{n,m})\Big)f(u_{n,m})v_{n,m}dx.
    \end{aligned}
\end{eqnarray*}

For $\delta > 0$ sufficiently small, we set
\begin{eqnarray*}
    \begin{aligned}\displaystyle
    \overline{u_{n,m}}=(Q\pi-\delta)^{\frac{1}{2}}\frac{u_{n,m}}{\|u_{n,m}\|_{1/2}},\ \ \ \overline{v_{n,m}}=(Q\pi-\delta)^{\frac{1}{2}}\frac{v_{n,m}}{\|v_{n,m}\|_{1/2}}.
    \end{aligned}
\end{eqnarray*}
Using Lemma $\ref{LEca}$ with $s=g(v_{n,m})/\sqrt{\beta_{0}}$ and $t=\overline{u_{n,m}}\sqrt{\beta_{0}}$, we have
\begin{eqnarray*}
    \begin{aligned}\displaystyle
    (Q\pi-\delta)^{\frac{1}{2}}\|u_{n,m}\|_{1/2}=&\int_{\mathbb{R}}\Big(I_{\mu_{1}}\ast G(v_{n,m})\Big)g(v_{n,m})\overline{u_{n,m}}dx\\
     \leq& \int_{\{x\in \mathbb{R}; \frac{g(v_{n,m})}{\sqrt{\beta_{0}}}\leq e^{1/4}\}} \Big(I_{\mu_{1}}\ast G(v_{n,m})\Big) (e^{\beta_{0}\overline{u_{n}}^{2}}-1)\\
    & + \int_{\{x\in \mathbb{R}; \frac{g(v_{n,m})}{\sqrt{\beta_{0}}}\geq e^{1/4}\}}\Big(I_{\mu_{1}}\ast G(v_{n,m})\Big) \frac{g(v_{n,m})}{\sqrt{\beta_{0}}}\Big(\log (\frac{g(v_{n,m})}{\sqrt{\beta_{0}}})\Big)^{\frac{1}{2}}\\
    & + \frac{1}{2} \int_{\{x\in \mathbb{R}; \frac{g(v_{n,m})}{\sqrt{\beta_{0}}}\leq e^{1/4}\}}\Big(I_{\mu_{1}}\ast G(v_{n,m})\Big) \frac{(g(v_{n,m}))^{2}}{\beta_{0}}\\
    :=& I_{1}+I_{2}+I_{3}.
    \end{aligned}
\end{eqnarray*}

Since $\|u_{n,m}\|^{2}=Q\pi-\delta \leq \frac{1+\mu_{1}}{2\beta_{0}}-\delta$, arguing as in Lemma $\ref{LEcb}$, we know that $I_{1},I_{3}\rightarrow 0$ as $n\rightarrow\infty$. By the exponential growth of $f$, given $\varepsilon>0$, there is a positive constant $C_{\varepsilon}$ such that $g(t)\leq C_{\varepsilon}e^{(\beta_{0}+\varepsilon)t^{2}}$, for all $t\geq 0$. Thus,
\begin{eqnarray*}
    \begin{aligned}\displaystyle
I_{2}&\leq \int_{\mathbb{R}}\Big(I_{\mu_{1}}\ast G(v_{n,m})\Big) \frac{g(v_{n,m})}{\sqrt{\beta_{0}}}\Big(\log (\frac{g(v_{n,m})}{\sqrt{\beta_{0}}})\Big)^{\frac{1}{2}}\\
&\leq \frac{1}{\sqrt{\beta_{0}}} \int_{\mathbb{R}}\Big(I_{\mu_{1}}\ast G(v_{n,m})\Big) g(v_{n,m})\Big(\log (\frac{C_{\varepsilon}}{\sqrt{\beta_{0}}}e^{(\beta_{0}+\varepsilon)v_{n,m}^{2}})\Big)^{\frac{1}{2}}\\
& \leq\frac{1}{\sqrt{\beta_{0}}} \int_{\mathbb{R}}\Big(I_{\mu_{1}}\ast G(v_{n,m})\Big) g(v_{n,m}) \Big(\log (\frac{C_{\varepsilon}}{\sqrt{\beta_{0}}}\Big)+(\beta_{0}+\varepsilon)v_{n,m}^{2}\Big)^{\frac{1}{2}}\\
& \leq o(1)+(1+\frac{\varepsilon}{\beta_{0}})^{1/2}\int_{\mathbb{R}}\Big(I_{\mu_{1}}\ast G(v_{n,m})\Big) g(v_{n,m})v_{n,m}.
    \end{aligned}
\end{eqnarray*}
Therefore, we obtain
\begin{eqnarray*}
    \begin{aligned}\displaystyle
(Q\pi-\delta)^{\frac{1}{2}}\|u_{n,m}\|_{1/2}\leq (1+\frac{\varepsilon}{\beta_{0}})^{\frac{1}{2}}\int_{\mathbb{R}}\Big(I_{\mu_{1}}\ast G(v_{n,m})\Big)g(v_{n,m})v_{n,m}dx+o(1).
    \end{aligned}
\end{eqnarray*}
Following the same steps as above, we are able to obtain that
\begin{eqnarray*}
    \begin{aligned}\displaystyle
    (Q\pi-\delta)^{\frac{1}{2}}\|v_{n,m}\|_{1/2}\leq (1+\frac{\varepsilon}{\alpha_{0}})^{\frac{1}{2}}\int_{\mathbb{R}}\Big(I_{\mu_{2}}\ast F(u_{n,m})\Big)f(u_{n,m})u_{n,m}dx+o(1).
    \end{aligned}
\end{eqnarray*}
In view of Lemma $\ref{LEcd}$, $u_{n}\rightarrow0$ strongly in $L^{p}(\mathbb{R})$, for $2<p<\infty$. Similarly to \cite{ACTT}, we may conclude that
\begin{eqnarray*}
    \begin{aligned}\displaystyle
    \int_{\mathbb{R}}\Big(I_{\mu_{2}}\ast F(u_{n,m})\Big)F(u_{n,m})dx\rightarrow 0 \ and \
    \int_{\mathbb{R}}\Big(I_{\mu_{1}}\ast G(v_{n,m})\Big)G(v_{n,m})dx\rightarrow 0,
    \end{aligned}
\end{eqnarray*}
when $(u_{n,m},v_{n,m})\rightharpoonup (0,0)$.

From $(\ref{PRCa})$, $c_{n,m}\leq Q\pi-\delta$ implies that
\begin{eqnarray*}
    \begin{aligned}\displaystyle
    |\int_{\mathbb{R}}((-\Delta)^{\frac{1}{4}} u_{n,m}(-\Delta)^{\frac{1}{4}} v_{n,m}+u_{n,m}v_{n,m})|\leq Q\pi-\delta +o(1).
    \end{aligned}
\end{eqnarray*}
Then
\begin{eqnarray*}
    \begin{aligned}\displaystyle
    (Q\pi-\delta)^{\frac{1}{2}}\|u_{n,m}\|_{1/2}+(Q\pi-\delta)^{\frac{1}{2}}\|v_{n,m}\|_{1/2} \leq& (1+\frac{\varepsilon}{\beta_{0}})^{\frac{1}{2}}\int_{\mathbb{R}}\Big(I_{\mu_{1}}\ast G(v_{n,m})\Big)g(v_{n,m})v_{n,m}dx\\
    &+(1+\frac{\varepsilon}{\alpha_{0}})^{\frac{1}{2}}\int_{\mathbb{R}}\Big(I_{\mu_{2}}\ast F(u_{n,m})\Big)f(u_{n,m})u_{n,m}dx\\
    \leq& 2\Big(1+\frac{\varepsilon}{\min\{\alpha_{0},\beta_{0}\}}\Big)^{\frac{1}{2}}(Q\pi-\delta)\\
    \leq& 2\Big(Q\pi-\frac{\delta}{2}\Big),
    \end{aligned}
\end{eqnarray*}
for $\varepsilon>0$ sufficiently small and $n$ sufficiently large.
Thus,
\begin{eqnarray*}
    \begin{aligned}\displaystyle
    \|u_{n,m}\|_{1/2}+\|v_{n,m}\|_{1/2}\leq2\Big(Q\pi-\frac{\delta}{2}\Big)^{\frac{1}{2}}.
    \end{aligned}
\end{eqnarray*}
Without loss of generality, we can assume that
\begin{eqnarray*}
    \begin{aligned}\displaystyle
    \|u_{n,m}\|_{1/2}\leq (Q\pi-\delta)^{\frac{1}{2}}\leq \Big(\frac{(1+\mu_{2})\pi}{2\alpha_{0}}-\frac{\delta}{2}\Big)^{\frac{1}{2}}.
    \end{aligned}
\end{eqnarray*}

Taking $t>1$ with $(H_{5})$ such that
\begin{equation}\label{PRCd}
    \begin{aligned}\displaystyle
    (\alpha_{0}+\varepsilon)(\frac{(1+\mu_{2})\pi}{2\alpha_{0}}-\frac{\delta}{2})(\frac{2}{1+\mu_{2}})t<\pi, \end{aligned}
\end{equation}
Then, using the Hardy-Littlewood-Sobolev and $H\ddot{o}lder$ inequalities, we get
\begin{eqnarray*}
    \begin{aligned}\displaystyle
    \int_{\mathbb{R}}\Big(I_{\mu_{2}}\ast F(u_{n,m})\Big)f(u_{n,m})u_{n,m}dx\leq C|F(u_{n,m})|_{\frac{2}{1+\mu_{2}}}|f(u_{n,m})u_{n,m}|_{\frac{2}{1+\mu_{2}}},
    \end{aligned}
\end{eqnarray*}
where $C$ is a positive constant. Then together with $(\ref{fnon})$, we obtain
\begin{eqnarray*}
    \begin{aligned}\displaystyle
    |f(u)|\leq\varepsilon |u|+b_{1}|u|^{q-1}(e^{(\alpha_{0}+\varepsilon) u^{2}}-1) \ \ \ \forall \ u \in H^{1/2}(\mathbb{R}).
    \end{aligned}
\end{eqnarray*}
Then,
\begin{eqnarray*}
    \begin{aligned}\displaystyle
    |f(u_{n,m})u_{n,m}|_{\frac{2}{1+\mu_{2}}}\leq\varepsilon |u_{n,m}|^{2}_{2}+C |u_{n,m}|^{\frac{1+\mu_{2}}{2t^{\prime}}}_{\frac{2qt^{\prime}}{1+\mu_{2}}}\Big(\int_{\mathbb{R}}[e^{\frac{2(\alpha_{0}+\varepsilon) t}{1+\mu_{2}}\|u_{n,m}\|^{2}(\frac{u_{n,m}^{2}}{\|u_{n,m}\|^{2}})}-1]dx\Big)^{\frac{1+\mu_{2}}{2t}},
    \end{aligned}
\end{eqnarray*}
where $t,t^{\prime}>1$ satisfying $\frac{1}{t}+\frac{1}{t^{\prime}}=1$. Together Trudinger-Moser inequality with $(\ref{PRCd})$, we deduce that
\begin{eqnarray*}
    \begin{aligned}\displaystyle
    \Big(\int_{\mathbb{R}}[e^{\frac{2(\alpha_{0}+\varepsilon) t}{1+\mu_{2}}\|u_{n,m}\|^{2}(\frac{u_{n,m}^{2}}{\|u_{n,m}\|^{2}})}-1]dx\Big)^{\frac{1+\mu_{2}}{2t}}\leq C_{1},
    \end{aligned}
\end{eqnarray*}
which implies that
\begin{eqnarray*}
    \begin{aligned}\displaystyle
    \int_{\mathbb{R}}\Big(I_{\mu_{2}}\ast F(u_{n,m})\Big)f(u_{n,m})u_{n,m}dx\leq \varepsilon^{2} |u_{n,m}|^{4}_{2}+ C_{2} |u_{n,m}|^{\frac{1+\mu_{2}}{t^{\prime}}}_{\frac{2qt^{\prime}}{1+\mu_{2}}}.
    \end{aligned}
\end{eqnarray*}

By $u_{n}\rightarrow0$ strongly in $L^{p}(\mathbb{R})$, for $2<p<\infty$, we get
\begin{eqnarray*}
    \begin{aligned}\displaystyle
    \int_{\mathbb{R}}\Big(I_{\mu_{2}}\ast F(u_{n,m})\Big)f(u_{n,m})u_{n,m}dx \rightarrow 0.
    \end{aligned}
\end{eqnarray*}
And then we conclude that $\|v_{n,m}\|\rightarrow 0$, repeating the same argument we get $\|u_{n,m}\|\rightarrow 0$. Hence, by case $1$, this situation cannot occur. Consequently, we obtain $(u_{0},v_{0})\neq 0$. So, $(u_{0},v_{0})$ is a nontrivial weak solution of $(\ref{a})$. Then choosing $\psi=u_{0}^{-}=\{-u_{0},0\}$ and $\varphi=0$ in $(\ref{PRCb})$, we obtain
\begin{eqnarray*}
    \begin{aligned}\displaystyle
    -\|u_{0}^{-}\|_{1/2}^{2}=\int_{\mathbb{R}}\Big(I_{\mu_{1}}\ast G(v_{0})\Big)g(v_{0})u_{0}^{-}\geq0,
    \end{aligned}
\end{eqnarray*}
which implies that $u_{0}^{-}=0$. Similarly, we deduce $v_{0}^{-}=0$, thus $u_{0},v_{0}\geq0$. If $u_{0}=0$, it is obvious to obtain $v_{0}=0$, allowing us to conclude that $u_{0},v_{0}>0$. Thus, the proof is complete.
\qed

\subsection*{Conflict of Interest}
The authors declared that they have no conflict of interest.

 \end{document}